\documentclass[ijoc,sglanonrev]{informs4}
\usepackage{eqndefns-left} % For checking the display equation width and equation environment definitions %
\RequirePackage{tgtermes}
\RequirePackage{newtxtext}
\RequirePackage{newtxmath}
\RequirePackage{bm}
\RequirePackage{endnotes}

\OneAndAHalfSpacedXI
\usepackage{multirow}
\usepackage{algpseudocode}
\usepackage{tikz}
\usepackage{enumitem}
\usepackage{epstopdf}
\usepackage{makecell}
\usepackage{comment}
\usepackage{amsfonts,amssymb,amsmath}
\usepackage{color,dsfont}
\usepackage{longtable}
\usepackage{rotating}
\usepackage{multicol}
\usepackage{lipsum}
\usepackage{bm}
\usepackage{tikz} 
\usepackage{xcolor}
\usepackage{hyperref}
\usepackage{graphicx}
\usepackage{caption}
\usepackage{booktabs}
\usepackage{threeparttable}

\usepackage{natbib}
\usepackage{amsmath}

\usepackage{accents}

\DeclareMathOperator{\conc}{conc}

\DeclareMathOperator{\conv}{conv}

\DeclareMathOperator{\Lovasz}{\texttt{L}}

\def\for{\mbox{ for }}

\def\MI{\text{MI}}

\def \x {\boldsymbol{x}}
\def \d {\boldsymbol{d}}
\def \y {\boldsymbol{y}}
\def \z {\boldsymbol{z}}
\def \s {\boldsymbol{s}}
\def \t {\boldsymbol{t}}

\def \R {\mathbb{R}}
\def \X {\mathcal{X}}
\def \one {\mathbf{1}}

\definecolor{links}{RGB}{204,36,29}

\usepackage{hyperref}
\hypersetup{
    colorlinks,
    citecolor=links,
    filecolor=links,
    linkcolor=links,
    urlcolor=links
}

\usepackage{mathtools}
\DeclarePairedDelimiterX{\infdivx}[2]{(}{)}{%
  #1\;\delimsize\|\;#2%
}
\newcommand{\infdiv}{D_{\text{KL}}\infdivx}

\usepackage{etoolbox}
\newenvironment{tightalign}
  {\begingroup
   \setlength{\abovedisplayskip}{0.6em}
   \setlength{\belowdisplayskip}{0.4em}
   \setlength{\abovedisplayshortskip}{0pt}
   \setlength{\belowdisplayshortskip}{0.6em}
  }
  {\endgroup\ignorespacesafterend}

\usepackage[normalem]{ulem}
\usepackage{natbib}
 \bibpunct[, ]{(}{)}{,}{a}{}{,}%
 \def\bibfont{\small}%
\EquationsNumberedThrough    % Default: (1), (2), ...
\TheoremsNumberedThrough     % Preferred (Theorem 1, Lemma 1, Theorem 2)
\ECRepeatTheorems  %  

\MANUSCRIPTNO{IJOC-0001-2024.00}

\begin{document}
%%%%%%%%%%%%%%%%

%\setcounter{page}{1}
%\DOI{10.1287/xxxx.0000.0000}%

% % Enter authors following the given pattern:
% \RUNAUTHOR{}
% \RUNAUTHOR{xxx}

% \RUNTITLE{Predictive Maintenance in Manufacturing}
\RUNTITLE{Optimal Information-Theoretic Data Reduction}
\RUNAUTHOR{He, Luo, and Zhao}

% \TITLE{Optimal Resource Allocation in Humanitarian Logistics: A Stochastic Programming Approach}
\TITLE{Optimal Data Reduction under Information-Theoretic Criteria}
%Globally Optimizing Data Reduction under Information Theoretic Criteria
%Optimal Data Reduction under Information Theoretic Criteria
\ARTICLEAUTHORS{%
%\AUTHOR{John Doe,\textsuperscript{a} Jane Smith,\textsuperscript{b}}
%\AFF{\textsuperscript{a}Department of Industrial Engineering, University of XYZ, \EMAIL{john.doe@xyz.edu; \textsuperscript{b}Department of Computer Science, University of ABC, \EMAIL{jane.smith@abc.edu}} 
\AUTHOR{Taotao He, Jun Luo, and Junkai Zhao}
\AFF{Antai College of Economics and Management,
Shanghai Jiao Tong University, \\ \EMAIL{hetaotao@sjtu.edu.cn} ( \href{http://taotaoohe.github.io}{\url{taotaoohe.github.io}}),  
\EMAIL{jluo\_ms@sjtu.edu.cn}, 
\EMAIL{zhaojunkai@sjtu.edu.cn}}

% \AUTHOR{Taotao He}
% \AFF{Antai College of Economics and Management,
% Shanghai Jiao Tong University, \EMAIL{hetaotao@sjtu.edu.cn}}

% \AUTHOR{Jun Luo}
% \AFF{Antai College of Economics and Management,
% Shanghai Jiao Tong University, \EMAIL{jluo_ms@sjtu.edu.cn}}
% \AUTHOR{Jun Luo}
% \AFF{Antai College of Economics and Management,
% Shanghai Jiao Tong University, \EMAIL{jluo_ms@sjtu.edu.cn}}

%
%\AUTHOR{Emily Johnson}
%\AFF{Department of Logistics,
%University of ABC, \EMAIL{emily.johnson@abc.edu}}
} % end of the block

\ABSTRACT{Selecting an optimal subset of features or instances under an information-theoretic criterion has become an effective preprocessing strategy for reducing data complexity while preserving essential information.  This study investigates two representative problems within this paradigm: feature selection based on the maximum-relevance minimum-redundancy criterion, and instance selection grounded in the Kullback–Leibler divergence. To address the intrinsic nonconvexities of these problems, we develop polyhedral relaxations that yield exact mixed-integer linear programming formulations, thereby enabling globally optimal data reduction. By leveraging modern optimization techniques, we further design efficient algorithmic implementations  capable of solving practically sized instances. Extensive numerical experiments on both real-world and synthetic datasets demonstrate that our method efficiently solves data reduction problems to global optimality, significantly outperforming existing benchmark approaches.
}%

% \FUNDING{This research was supported by [grant number, funding agency].}

%Supplemental Material:
%Data Ethics & Reproducibility Note:

%\KEYWORDS{Stochastic programming, Decision support,Uncertainty, Disaster response, Optimization}

%\HISTORY{Received: Month DD, YYYY; Accepted: Month DD, YYYY; Published Online: Month DD, YYYY}

\maketitle
%%%%%%%%%%%%%%%%%%%%%%%%%%%%%%%%%%%%%%%%%%%%%%%%%%%%%%%%%%%%%%%%%%%%%%

% Text of your paper here

\section{Introduction}

% \ins{Subset selection}

% Reducing the size of inputs of ML/DS algorithms is crucial. What do we mean by size? It incudlues feature, number of observations. 
% preprocessing data to reduce its complexity while retaining its essential information  plays a critical role in boosting training efficiency.

In recent years, the explosive growth of data has fueled remarkable advances in the fields of machine learning and data mining by enabling the training of increasingly sophisticated models and the discovery of complex patterns~\citep{wei2015submodularity, garcia2015data}. As data continues to be collected at an unprecedented pace, data reduction has emerged as a crucial preprocessing technique to reduce data complexity while retaining its essential information, thereby boosting the efficiency  of model training and knowledge discovery~\citep{zha2025data}. This strategy is often achieved by either reducing the number of features (i.e., feature selection), or reducing the number of instances (i.e., instance selection). Generally speaking, feature selection aims to identify a subset of the most relevant features from all available features in the dataset, which can facilitate data visualization, reduce storage and computation costs, and improve predictive performance~\citep{guyon2003introduction,li2017feature,zha2025data}. Instance selection, on the other hand, aims to retain a representative subset of the data, helping to alleviate computational constraints while maintaining learning quality~\citep{zha2025data}, and also grasp a better understanding of the whole dataset~\citep{zhang2021review}.

One of the challenges in data reduction lies in evaluating the quality of the selected features or instances. A natural approach is to leverage information-theoretic criteria to quantify how relevant or representative a subset is~\citep[see][and references therein]{brown2012conditional,li2017feature,wei2015submodularity,zhang2021review}. In the feature selection, many information-theoretic criteria are proposed to balance feature relevance and redundancy. A prominent example is the widely used Minimum-Redundancy Maximum-Relevance (mRMR) criterion proposed by~\cite{peng2005feature}, in which the average mutual information between each selected feature and the target variable quantifies relevance, while the average mutual information among pairs of selected features is used as a redundancy penalty. For instance selection, information-theoretic criteria such as the Kullback–Leibler (KL) divergence  are often used to assess the discrepancy between the distribution over all features of the original data and that of a selected data subset~\citep{wei2015submodularity,zhang2021review}.

Selecting an optimal subset of features or instances under a given information-theoretic criterion is often formulated as a nonlinear integer programming problem. Such formulation typically uses binary decision variables and linear constraints to represent all possible subsets, while the objective function—derived from the chosen information-theoretic measure—is inherently nonlinear. For instance, in the mRMR criterion, the objective includes bilinear terms representing mutual information between pairs of selected features, as well as fractional terms to compute average mutual information, as formally presented in~\eqref{eq:mrmr-fp}. In the case of instance selection based on KL divergence, the objective involves a composition of logarithmic and rational functions. This captures the discrepancy between the distributions of the original data and the selected subset, where the main component can be formally written as~\eqref{eq:logratio}.  These nonlinear structures introduced by information-theoretic measures significantly increase the complexity of the underlying subset selection problem.

Most existing approaches to solving information-theoretic data reduction problems yield near-optimal solutions~\citep{li2017feature,zha2025data}. For instance,~\cite{peng2005feature} introduce the mRMR criterion along with a forward greedy algorithm for incremental feature selection, while~\cite{brown2012conditional} later propose a backward greedy variant for the same class of problems. A different line of work by~\cite{naghibi2014semidefinite} employs a semidefinite programming relaxation to design a rounding algorithm that achieves near-optimality. Similarly,~\cite{nguyen2014effective} develop a rounding algorithm based on spectral relaxation. In the context of selecting a representative subset from large datasets for classifier training,~\cite{wei2015submodularity} demonstrate that the performance loss can be expressed as the difference between two submodular functions, allowing for approximate minimization via algorithms for difference-of-submodular-function optimization~\citep{iyer2012-min-diff-sub-funcs-arxiv,el2023difference}.  Recently,~\cite{zhang2021review} propose two approximation algorithms to select a subset of reviews that closely preserves the distribution of the original corpus in terms of KL divergence.

% adopt an incremental, greedy strategy, selecting one feature or sample at a time~\citep{li2017feature,zha2025data}. ~\cite{peng2005feature} propose the mRMR criterion and provide an incremental greedy algorithm for feature selection. Later, various \cite{brown2012conditional} develop incremental and backward greedy algorithms within a unifying framework for information-theoretic feature selection

While computationally efficient, these methods often result in suboptimal solutions due to their inherently greedy or myopic design.  Better performance gains can potentially be achieved by making global selection decisions~\citep{naghibi2014semidefinite}. In this paper, we focus on developing global optimization approaches for two representative data reduction problems: feature selection under the mRMR criterion~\citep{peng2005feature} and instance selection under KL divergence~\citep{zhang2021review}. Specifically, we propose mixed-integer linear programming (MIP) formulations for both problems.  Our formulations are built on convex relaxations of the nonlinear structures commonly appearing in information-theoretic data reduction problems.  Thus, our methodology also applies naturally to other settings, such as the feature selection under correlation-based measures~\citep{yu2003feature} and the instance selection  considered in~\cite{wei2015submodularity}.  By leveraging modern MIP solvers, such as \texttt{Gurobi}, and efficient cutting-plane algorithms, our approach can globally solve practically sized instances within a reasonable computational time. 

% can be applied to solve subsection selections problems, e.g., joint mutual information (JMI) \citep{yang1999data}, correlation-feature-selection measure (CFS) \citep{yu2003feature} and redundancy rate (RED) \citep{zhao2010efficient}. 

\subsection{Contributions}

We highlight the contributions of our paper as follows:
%\ins{Our solution}
\begin{enumerate}
    \item 
    % We construct a polyhedral relaxation to tackle the two  nonconvex structures,  ratio and bilinearity, in feature selection under mRMR. Our relaxation yields a ``bound-free" formulation, while existing techniques produces MIP formulations using derived bounds on auxiliary variables~\citep{mehmanchi2021solving}. We also theoretically analyze the formulation and, in Theorem~\ref{theorem:comparison}, we show that the continuous relaxation of our formulation is tighter than that of an MIP formulation obtained using recursive McCormick envelopes considered in~\cite{mehmanchi2021solving}. 
We develop a polyhedral relaxation to address the two key sources of nonconvexity—ratios and bilinearity—in feature selection under the mRMR criterion. Our relaxation, together with binary variables for modeling selected features,  yields an MIP formulation.
% Unlike existing approaches that construct MIP formulations by introducing bounds on auxiliary variables~\citep{mehmanchi2021solving}, our method results in a ``bound-free" formulation. 
We also provide a theoretical analysis and prove in Theorem~\ref{them:mRMR-Pers} that the continuous relaxation of our formulation is tighter than the MIP model based on recursive McCormick envelopes that appeared in~\cite{mehmanchi2021solving}.

    \item  
    % In Section~\ref{section:irssformulation}, we propose a polyhedral relaxation for the log-rational function in the review selection under KL divergence, proposed in~\cite{zhang2021review}. Our relaxation exactly represents the log-rational function over all feasible subsets, thus yielding an MIP formulation for the review selection problem (see Theorem~\ref{theorem:envelope}). To the best of our knowledge, exact MIP formulations for the log-rational function have not been previously proposed in either machine learning or the operations research community. 
    We present a polyhedral relaxation for the log-rational function, defined as in~\eqref{eq:logratio}, arising in the instance selection problem under KL divergence introduced by~\cite{zhang2021review}. Our relaxation exactly represents the log-rational function, leading to an exact MIP formulation for the problem (see Theorem~\ref{theorem:envelope}). To the best of our knowledge, exact MIP formulations for the log-rational function have not been previously proposed in either machine learning or operations research communities.

    \item  
% We perform a comprehensive computational evaluation of our formulations for globally solving data reduction problems. For feature selection under the mRMR criterion, we use ten real-world datasets with feature counts ranging from 19 to 856, where mRMR is widely regarded as one of the most effective selection criteria (see Table~\ref{tab: dataset}). Our bound-free formulation achieves provable optimality in substantially less time and yields significantly tighter relaxations compared to the three alternative formulations proposed in~\cite{mehmanchi2021solving}. Specifically, our formulation solves small- to medium-sized instances in under 10 seconds and medium- to large-sized instances in under 100 seconds, while all three alternatives fail to solve even small instances within 3600 seconds (Table~\ref{tab: mRMR_real_data}). This computational advantage is further corroborated by experiments on synthetic datasets (Table~\ref{tab: mRMR_synthetic}).    
    We conduct a  comprehensive computational evaluation of our formulations for globally solving data reduction problems. For feature selection under  the mRMR criterion, we use ten real-world datasets ranging from 19 to 856 features, where mRMR is recognized as one of the most effective selection criteria (see Table~\ref{tab: dataset} and Appendix~\ref{apx: mRMR_real_data}).    
    % We computationally evaluate the performance of our formulations on globally solving data reduction problems. For feature selection under mRMR, we collect ten real datasets, ranging from 19 features to 856 features, on which mRMR is shown to be one of the best criteria for feature selection  (see Table~\ref{tab: dataset}).
    % global optimal solutions under mRMR yield performance improvements (see Table~\ref{tab: dataset}). 
    Our formulation  achieves provable optimality in substantially less time and yields significantly tighter relaxations compared to the three alternative formulations evaluated in~\cite{mehmanchi2021solving}. Specifically, our formulation solves small- to medium-sized instances within ten seconds, and medium- to large-sized instances within a few hundred seconds, while the three alternatives fail to solve even small-sized instances within 3600 seconds (Table~\ref{tab: mRMR_real_data}). This computational advantage is further confirmed by experiments on synthetic datasets (Table~\ref{tab: mRMR_synthetic}).
% For the \replace{review selection}{instance selection} problem, our formulation, together with a cutting-plane implementation, exactly solves small- to medium-sized instances within a few seconds. In contrast, using the same amount of time, local search heuristics proposed in~\cite{zhang2021review} produce average final optimality gaps of approximately 20\%. Furthermore, our approach can solve larger instances, enabling global optimization at scales previously unattainable.
For the instance selection problem, our formulation, together with a cutting-plane implementation, solves small- to medium-sized instances exactly within a few seconds. In contrast, within the same time budget, the local search heuristics proposed in~\cite{zhang2021review} yield average final optimality gaps of above 17\%. Additionally, our approach scales to larger instances, enabling global optimization at levels previously considered intractable.
    % For the review selection problem, our formulation, together with a cutting-plane implementation, exactly solves small to middle sized instances within few seconds, while, within a similar amount of computational time, the average final gaps produced local search heuristics proposed in~\cite{zhang2021review} are near 20\%. Moreover, our approach allows for global solutions of larger sized instance.     
\end{enumerate}
Our paper demonstrates how convex relaxation techniques can be leveraged to preprocess data, reducing data complexity while preserving essential information. This approach aligns with the growing body of research applying modern mixed-integer linear and nonlinear programming to machine learning tasks~\citep{bertsimas2019machine, huchette2023deep, tillmann2024cardinality}. Notably, convex relaxations have proven effective in training various machine learning models, including sparse regression~\citep{bertsimas2020sparse, gomez2021mixed, atamturk2025rank}, sparse principal component analysis~\citep{d2007direct, bertsimas2022solving, dey2022using, kim2022convexification, li2025exact}, sparse and low rank matrix decomposition~\citep{bertsimas2023sparse} and experiment design~\citep{li2025strong}, and in optimizing trained machine learning models such as neural networks~\citep{anderson2020strong,kronqvist2025p}, decision tree~\citep{mivsic2020optimization,kim2024reciprocity}, and optimization with learned constraints~\citep{maragno2025mixed}.

% and decision tree~\citep{mivsic2020optimization, mistry2021mixed}.

% For instance, \cite{atamturk2025rank} introduce a new semidefinite programming (SDP) relaxation for sparse regression that is tighter than the perspective relaxation of \cite{bertsimas2020sparse}, while \cite{gomez2021mixed} present a convexification approach for sparse regression under various information-theoretic criteria. Since the introduction of an SDP relaxation for sparse principal component analysis by \cite{d2007direct}, several stronger relaxations have been developed \citep{bertsimas2022solving, dey2022using, kim2022convexification, li2025exact}.~\cite{bertsimas2023sparse} propose a novel SDP relaxation for sparse plus low rank matrix decomposition which finds various applications in machine learning.~\cite{khajavirad2025inference} propose various polyhedral relaxations for inference in higher-order undirected graphical models with binary labels.~\cite{li2025exact} leverage convex envelope techniques to solve the regularized A-optimal design \cmt{Jun}{is that well known?}. \cmt{Jun}{from my point of view, I am not so sure what is the purpose of this paragraph. Of course, I know it demonstrate our method has important role in ML. But seems the following reference is too heavy. But I don't know how to better place it.}
\subsection{Structure}
The rest of the paper is organized as follows. In Section~\ref{sec: subset_selection}, we formally define the feature reduction under mRMR and the instance reduction under KL divergence.  In Section~\ref{sec: IP_formulation}, we present our MIP formulations for both problems. In Section~\ref{section:implementation}, we discuss the implementations of our formulations and present computational experiments in Section~\ref{sec: numerical}. Finally, conclusions are given in Section~\ref{section:conlcusion}. Supplementary experimental details and proofs are provided in the E-Companion.

\section{Models and Preliminaries}
\label{sec: subset_selection}
Before formally introducing data reduction problems studied in this paper, we briefly review key information-theoretic concepts that are useful for quantifying the quality of the selected subsets. For a more comprehensive treatment, see~\cite{polyanskiy2025information}. Let $X$ be a discrete random variable with probability mass function (PMF) $p(\cdot)$ over a finite set  $\mathcal{X}$. The entropy of $X$, which measures its intrinsic uncertainty, is defined as
\begin{equation*}
H(X) := \sum_{x \in \mathcal{X}} p(x)\log \frac{1}{p(x)}.
\end{equation*}
The conditional entropy of $X$ given another discrete random variable $Y$, with joint PMF $p(x,y)$, is 
\begin{equation*}
    H(X \mid Y) :=  \sum_{x\in \mathcal{X}, y\in \mathcal{Y}} p(x, y)\log  \frac{1}{p(x\mid y)},
\end{equation*}
where  $p(\cdot \mid y)$ is the conditional PMF of $X$ given $Y = y$. The mutual information between $X$ and $Y$ is used to measure the amount of information shared by $X$ and $Y$ and is defined as :
\begin{equation*}
    \MI(X,Y) := H(X) - H(X\mid Y) = \sum_{x\in \mathcal{X}, y\in \mathcal{Y}}p(x, y)\log\left[\frac{p(x, y)}{p(x)p(y)}\right].
\end{equation*}
All these quantities can be estimated from data; we refer readers to~\cite{paninski2003estimation} and Section 3.3 of~\cite{brown2012conditional} for estimation methods.
% To use these in data reduction problems, we should be able to estimate them from the data. For this we refer the readers to~\cite{paninski2003estimation} and Section 3.3~\cite{brown2012conditional}. 
For a pair of PMFs $p(\cdot)$ and $q(\cdot)$ with a common support $\mathcal{X}$, the KL divergence between $p$ and $q$ is:
\begin{equation*}
     \infdiv{p}{q} = \sum_{x\in \mathcal{X}}p(x)\log\left[\frac{p(x)}{q(x)}\right],
\end{equation*}
which quantifies the extent to which $p$ diverges from the distribution $q$. 

\subsection{Feature Selection under mRMR}
\label{sec: feature_subset_selection}
% Recent survey~\cite{li2017feature}.
% We start with introducing some elementary concepts from information theory. The entropy of a discrete random variable $X$ is defined as $H(X) := -\sum_{x\in \mathcal{X}} p(x)\log p(x)$,
% where $p(\cdot)$ denotes the probability distribution function and $\mathcal{X}$ denotes the support set of $X$. It quantifies the average level of uncertainty or information associated with the variable $X$. The conditional entropy for a discrete random variable $X$ given another variable $Y$ is
% \begin{equation*}
%     H(X \mid Y) := -\sum_{x\in \mathcal{X}, y\in \mathcal{Y}} p(x, y)\log p(x\mid y),
% \end{equation*}
% where $p(\cdot)$ and $p(\cdot\vert \cdot)$ denote the probability distribution function  and the conditional probability function, and $\mathcal{Y}$ is the support set of $Y$. It measures the amount of information left in $X$ based on having knowledge regarding the random variable $Y$.  Therefore, the information shares by variables $X$ and $Y$ can be measured by mutual information:
% \begin{equation*}
%     \MI(X,Y) = H(X) - H(X\mid Y) = \sum_{x\in \mathcal{X}, y\in \mathcal{Y}}p(x, y)\log\left[\frac{p(x, y)}{p(x)p(y)}\right],
% \end{equation*}
% which reflects the dependence or relevance of $X$ and $Y$. 

Suppose that there are $m$ features in a dataset. Let $[m]:=\{1, 2, \ldots, m\}$, and let $\{\gamma_j\}_{j \in [m]}$ be the set of $m$ features, where $\gamma_j$ represents the $j^{\text{th}}$ feature. Given the target variable $Y$ to be predicted in a supervised learning task, $\MI(\gamma_j,Y)$ quantifies the relevance of the $j^{\text{th}}$ feature to $Y$, and $\MI(\gamma_j,\gamma_k)$ measures the redundancy between two features $\gamma_j$ and $\gamma_k$. Using these notions, \cite{peng2005feature} introduce the mRMR criterion, which leads to the following feature selection problem:
% the feature selection problem under the Minimum-Redundancy Maximum-Relevance (mRMR) criterion, proposed by~\cite{peng2005feature}, is defined as follows:
\begin{tightalign}
\begin{align}\label{eqn: mRMR}
\max_S \biggl\{I_{\text{mRMR}}(S):=
\frac{1}{\vert S \vert}\sum_{j\in S} \MI(\gamma_j, Y) - \frac{1}{\vert S\vert^2}\sum_{j, k \in S} \MI(\gamma_j, \gamma_k) \biggm| L\leq \vert S \vert \leq U,\ S \subseteq [m] \biggr\},   \tag{\textsc{mRMR}}
\end{align}
\end{tightalign}
% \noindent
where $\vert \cdot \vert$  denotes the cardinality of a subset,  and $L$  and $U$ are positive integers ensuring that enough features are retained to preserve predictive power while balancing computational cost. The first term in the objective $I_{\text{mRMR}}(S)$ captures the average relevance of the feature subset $S$ to $Y$, i.e., the average information from $S$ that helps explain $Y$, and the second term penalizes  the average redundancy among features in $S$.

\subsection{Instance Selection with KL Divergence}
\label{sec: data_subset_selection}

% Data subset selection refers to the process of finding a small subset of data which represents the characteristics of the original full data.
% Let $\mathcal{D} = \{d_j, j = 1,\ldots, m\}$ be the full data, where $n$ denotes the total number of observations, information-theoretic data subset selection aims to select the most informative $
% \mathcal{S}\in \mathcal{D}$ according to some information criteria.

For instance selection, we consider the online reviews selection setting studied in~\cite{zhang2021review}. 
%study a modified KL divergence based criterion proposed in~\cite{zhang2021review},  applied within the framework of the review selection problem described therein. 
Formally, we represent the full review dataset as a binary matrix $D_{n\times m}$, where each row $i \in [n]$ corresponds to a review (an instance) and each column $j \in [m]$ corresponds to an opinion (a feature). Here, $n$ denotes the total number of reviews, and $m$ denotes the total number of distinct opinions. Each entry $d^i_j$ in $D$ is defined such that  $d^i_j = 1$ if the $i^{\text{th}}$ review contains the $j^{\text{th}}$ opinion, and $d^i_j = 0$ otherwise.
The objective of instance/review selection in this setting is to select a representative subset of reviews $S \subseteq [n]$ that could cover as
many opinions as possible in the review corpus (i.e., full review data $D$), with its distribution of all opinions being largely consistent with that of the corpus.

% Instance selection attempts to select a subset of instances $S \subseteq [n]$ that retain the original feature characteristics of $D$. 

% For instance selection, we study a modified KL divergence based criteria proposed in~\cite{zhang2021review},  applied within the framework of the review selection problem described therein. Formally, we represent the full dataset using a matrix $D_{n\times m}$, where each row $i\in [n]$ represents a single instance while each column $j\in [m]$ stands for a feature, with
% $n$ (resp. $m$) denoting the number of instances  (resp. features). Instance selection attempts to select a
% subset of instances $S \subseteq [n]$ that retain the original feature characteristics of $D$.

% In the context of review selection, an instance corresponds to an individual review and a feature corresponds to an opinion reflected in the review.
% Each entry $d^i_j$ in $D$ is defined such that  $d^i_j = 1$ if the $i^{\text{th}}$ review contains the $j^{\text{th}}$ opinion, and $d^i_j = 0$ otherwise.
% The objective of review selection is to select a representative
% subset of reviews $S \subseteq [n]$ that could cover as
% many opinions as possible in the review corpus (i.e., full review data $D$), with its opinion distribution over all of the features being largely consistent with that of the corpus.

To measure the representativeness of the selected review subset, let $P^S_j$ denote the proportion of a given opinion $j$ that occurs in the review subset $S$, which is
\begin{equation*}
    P_{j}^{S} := \frac{1}{\vert S\vert}\sum_{i \in S}d^i_{j},
\end{equation*}
and, for simplicity, let  $P_j$  denote $P^{[n]}_j$. Without loss of generality, we assume that for each opinion $j \in [m]$, its score $P_j$  is not zero since otherwise we do not consider the $j^{\text{th}}$ opinion in our model. Then, the review subset selection with modified KL divergence (RSKL) is defined as
\begin{tightalign}
\begin{align}\label{eqn: Entropy_review_selection}
\min_S \biggl\{ I_{\text{KL}}(S):= \sum_{j\in [m]}P_{j}\Big\vert\log\frac{P_{j}}{P^S_{j}}\Big\vert \biggm| L\leq \vert S \vert \leq U \biggr\},   
 \tag{\textsc{RSKL}}
 \end{align}
\end{tightalign}
where the objective function is the aggregated absolute KL divergence between the distribution of the $j^{\text{th}}$ opinion in the subset $S$ and that in the full review data.
To ensure correctness, if $P_{j}^S = 0$ and $P_{j} > 0$, then 
\begin{equation}\label{eq:penalty}
P_j\Big\vert \log\frac{P_{j}}{P_{j}^S} \Big\vert = \delta := \max\left\{m+1, n\right\} \cdot \log (n),
\end{equation}
where $\delta$ represents a large penalty for the exclusion of opinion $j$ from $S$. 
% Let $\vert S\vert = K < n$, we need to guarantee that $\delta > \sum_j P_j\Big\vert\log(K P_j)\Big \vert - \sum_j P_j\Big\vert\log((K+1)P_j/2)\Big \vert$?
 % The main difference between the modified relative entropy in the objective of Problem~(\ref{eqn: Entropy_review_selection}) and relative entropy in Equation~(\ref{eqn: relative_entropy}) is the absolute value operator outside the log function. 
 In~(\ref{eqn: Entropy_review_selection}), if $P_{j}^S$  significantly deviates from $P_{j}$, the term $\Big\vert\log\frac{P_{j}}{P_{j}^S}\Big \vert$ will move far away from 0, leading to a substantially high objective value.  \cite{zhang2021review} examine and show that  the reviews selected by~(\ref{eqn: Entropy_review_selection}) perform favorably across various evaluation metrics and user feedback, which validates that the proposed formulation effectively addresses user needs for the selection of informative reviews.

% Solving Problem~\eqref{eqn: Entropy_review_selection} to globally optimality is also challenging due to the multiple nonlinear and nonsmooth items as well as the combinatorial nature. Similar structure can be found in many information theoretic subset selection problems, such as data selection for Na\"ive Bayes \citep{wei2015submodularity} and causal inference \citep{sun2016mutual}. \cite{zhang2021review} prove that problem~(\ref{eqn: Entropy_review_selection}) is NP-hard and propose a heuristic method called combined search (ComS) to solve it. ComS treats the problem as two optimization steps, the first of which solves the continuous version of the problem by relaxing the integer constraints and the second of which transforms the continuous solution into a discrete one. However, as a heuristic method, ComS cannot guarantee the global optimality of the selected subset.

%\subsection{Other Related Literature}

% It has many applications in data science, 
% predictive performance of a classifier trained on it is close to that of a classifier trained on the full training data.

% \begin{itemize}
% \item Mutual information in causal inference \citep{sun2016mutual}

% \begin{align}
%         \max_{x} \quad & \log(\sum_i x_i+\vert T\vert) + \frac{1}{\sum_i x_i + \vert T\vert}\sum_{b\in B}\Big(\sum_i m_{ib}x_i\Big[\log(\sum_i m_{ib}x_i) \\
%         &- \log(T_b + \sum_i m_{ib}x_i) - \log(\sum_i x_i)\Big]-T_b\log(T_b+\sum_i m_{ib}x_i)\Big)\\
%     \text{s.t.} \quad &\sum_{i\in N}x_i \le K,  
%     x \in \{0, 1\}^n    
%  \label{eqn: MI_causal_analysis}
%  \end{align}

% \end{itemize}

\section{Integer Programming Formulations via Convex Relaxations}
\label{sec: IP_formulation}
% In this section, we present our MIP formulations for data reduction problems,~\eqref{eqn: mRMR} and~\eqref{eqn: Entropy_review_selection}, introduced in Section~\ref{sec: subset_selection}. Our formulations are built upon polyhedral relaxations for nonlinear information-theory based criteria.  

In this section, we develop MIP formulations for the data reduction problems~\eqref{eqn: mRMR} and~\eqref{eqn: Entropy_review_selection} introduced in Section~\ref{sec: subset_selection}. Our approach employs polyhedral relaxations to handle the inherent nonconvexity of the  information-theoretic objectives. Specifically, in Section~\ref{sec: mRMR_formulation}, we formulate~\eqref{eqn: mRMR} as a fractional optimization problem, allowing us to exploit recent advances in fractional programming~\citep{he2024convexification}. In Section~\ref{section:irssformulation}, we express the nonconvex objective in~\eqref{eqn: Entropy_review_selection} as the difference of two composite functions, for which convex and concave envelopes can be derived.

% Specifically, in section~\ref{sec: mRMR_formulation}, we leverage recent convexification techniques in fractional programming to derive 
% In section~\ref{sec: mRMR_formulation}, we describe our MIP model for mRMR optimization problem and theoretically compare our formulation with existing alternate formulations.

\subsection{Feature Subset Selection under mRMR}\label{sec: mRMR_formulation}
% Let $\x = (x_1, x_2, \ldots, x_m) \in \{0,1\}^m$ be a vector of binary decision variables modeling the subset of features selected, i.e., $x_i = 1$ if and only if the feature $i$ is selected. Then, it follows readily that  the feature subset selection model~(\ref{eqn: mRMR}) can be expressed as the following binary fractional program:
% \begin{align}\label{eq:mrmr-fp}
%     \max_{\x} \Biggl\{  \frac{\sum_{i \in [m]} \sum_{j \in [m]} \bigl(\MI(f_i, Y) - \MI(f_i, f_j)\bigr) \cdot x_i x_j}{\sum_{i \in [m]} \sum_{j \in [m]} x_i x_j} \Biggm| L \le \sum_{i \in [m]}x_i \le U,\ \x \in \{0, 1\}^m  \Biggr\} .
%     % \text{s.t.} \quad &L \le \sum_{i \in [m]}x_i \le K,\\ 
%     % \quad & \x \in \{0, 1\}^m, 
% \tag{\textsc{mRMR-Frac}}
% \end{align} 
% Note that this formulation involves two types of nonlinear structures, namely, a rational function and  bilinear functions in both denominator and enumerator. Following convexification techniques in fractional programming~\citep{he2024convexification}, we use McCormick envelope and the reformulation-linearization technique (RLT) to handle bilinearity, and then use a perspective transformation to treat the rational structure.
% In Section~\ref{section:pers}, we present our MIP model, and in Section~\ref{section:comparison}, we theoretically investigate the strength or our formulation. 

Let $\x = (x_1, x_2, \ldots, x_m) \in \{0,1\}^m$ denote the vector of binary decision variables, where $x_i = 1$ if and only if the $i^{\text{th}}$ feature is selected. Using this notation, the feature subset selection model~\eqref{eqn: mRMR} can be equivalently reformulated as the following binary fractional program:
\begin{tightalign}
    \begin{align}\label{eq:mrmr-fp}
    \max_{\x} \Biggl\{  \frac{\sum_{i \in [m]} \sum_{j \in [m]} \bigl(\MI(\gamma_i, Y) - \MI(\gamma_i, \gamma_j)\bigr) \cdot x_i x_j}{\sum_{i \in [m]} \sum_{j \in [m]} x_i x_j} \Biggm| L \le \sum_{i \in [m]}x_i \le U,\ \x \in \{0, 1\}^m  \Biggr\},     % \text{s.t.} \quad &L \le \sum_{i \in [m]}x_i \le K,\\ 
    % \quad & \x \in \{0, 1\}^m, 
\tag{\textsc{mRMR-Frac}}
\end{align} 
\end{tightalign}
where the feasible region will be denoted as $\X$. Note that the objective function contains two sources of nonconvexity: (i) a fractional structure, and (ii) bilinear terms in both the numerator and denominator.  To address these challenges, we use recent convexification techniques from~\cite{he2024convexification} to construct a polyhedral relaxation for the objective function. This relaxation, together with binary variables for modeling selected features, yields an exact MIP formulation. 

% Specifically, the bilinear terms are handled using McCormick envelopes~\citep{mccormick1976computability}, while the rational structure is treated through the perspective transformation in~\cite{he2024convexification}. 

Before presenting our formation, we discuss prevalent MIP formulations in the literature. To linearize the objective function of ~\eqref{eq:mrmr-fp},  it is often to introduce the following auxiliary variables:
\begin{equation}\label{eq:transformation}
    \rho = \frac{1}{\sum_{s,t \in [m]} x_sx_t},   \quad y_i = \frac{x_i}{\sum_{s,t\in [m]}  x_sx_t} \for i \in [m], \quad \text{ and }  \quad z_{ij} = \frac{x_ix_j}{\sum_{s,t \in [m]} x_sx_t} \for i, j \in [m].
\end{equation}
With these definitions,~\eqref{eq:mrmr-fp} is equivalent to the following mixed-binary trilinear programming:
\begin{equation}\label{mrmr-tri}
\begin{aligned}
     \max \Biggl\{ \sum_{i, j \in [m]}  (\MI(\gamma_i, Y) - \MI(\gamma_i, \gamma_j) )  z_{ij} \Biggm| \x \in \X,\ \sum_{i,j \in [m]}z_{ij} = 1,\   z_{ij} = x_ix_j\rho \quad \forall i, j \in [m]\Biggr\}.
\end{aligned}
\end{equation}
Now, linearizing the resulting cubic terms, which  involve the products of two binary and one continuous variables, yields MIP formulations for~\eqref{eqn: mRMR}. In particular, recursively using McCormick envelopes~\citep{mccormick1976computability} to  relax   $y_i = x_i\rho$ and then $z_{ij} = y_ix_j$ leads to the following model: 
% where $\rho^L = \frac{1}{K^2}$ and $\rho^U = \frac{1}{L^2}$. 
% where $y_i^L = 0$ and $y_i^U = \frac{1}{L^2}$. 
\begin{tightalign}
\begin{align}\label{eq:mrmr-RMC}
   \begin{aligned}
        \max \quad & \sum_{i, j} \bigl(\MI(\gamma_i, Y) - \MI(\gamma_i, \gamma_j)\bigr) \cdot z_{ij} \notag \\
    \text{s.t.} 
        \quad &\x \in \X ,\ \sum_{i,j \in [m]}z_{ij} = 1,\ z_{ii} =y_i\quad \text{ for } i \in [m] \notag  \\
   &     \max\bigl\{\rho^Lx_i, \rho^Ux_i + \rho -\rho^U \bigr\} \leq y_i \leq \min \bigl\{\rho^L x_i + \rho - \rho^L, \rho^Ux_i \bigr\} \quad \for i \in [m]
% \label{eq:mRMR-RMC-1} 
\nonumber\\
   &  \max \bigl\{0, \rho^Ux_j + y_i -\rho^U \bigr\} \leq z_{ij} \leq \min \bigl\{  y_i , \rho^Ux_j \bigr\}  \quad \for i \in [m] \text{ and } j \in [m]\setminus \{i\}, 
        % \quad &   w_{iik}  = z_{ik} && \for k, i \in [m] \label{mRMR_full_cons-5} \\
            % \quad & 0  \le w_{ijk} \le y_k \text{ and } y_k + z_{ij} - \rho \le w_{ijk} \le z_{ij} && \for i \neq j,k \in [m] \label{mRMR_full_cons-6} \\
           % \quad & \rho (\sum_i a_i x_i + \sum_{i}\sum_{j} b_{ij}x_i x_j) \ge 1\\
    % \quad & \sum_{i} w_{ij} - y_j \ge 0 \for j \in [m]
    % \\
    % \quad & \sum_{i} y_{i} - \rho - \sum_{i,j}w_{ij} +y_j \ge 0 \for j \in [m]\\
    % \quad & \rho \ge 0 \quad 0 \leq y \leq \rho \text{ and } \quad \rho L \leq  \sum_{i \in [M]}y_i \leq \rho K \label{mRMR_full_cons-8}
 \end{aligned} \tag{\textsc{mRMR-Rmc}}
\end{align}
\end{tightalign}
where $\rho^L$ and $\rho^U$  are two constants such that $\rho^L \leq \rho \leq \rho^U$ for every $\x \in \X$. This model has been studied in~\cite{mehmanchi2021solving}, and other alternative models in~\cite{mehmanchi2021solving} are discussed in Appendices~\ref{apx: mRMR-bigm} and~\ref{apx: mRMR-vd}. In the following, we present our formulation and theoretically show that our formulation is tighter than~\eqref{eq:mrmr-RMC}. 

% The computational performance of~\eqref{eq:mrmr-RMC} depends on the tightness of bounds on $\rho$. Instead of treating $z_{ij}$ as 

% In Section~\ref{section:computation_mrmr}, we conduct computational experiments to show that our formulation  significantly outperforms existing ones. Here, we focus on presenting a theoretical comparison between our model and~\eqref{eq:mrmr-RMC}.  

\subsubsection{Perspective reformulations.}\label{section:pers} 

Note that $\rho > 0$ holds for all feasible $\x$. Following Theorem 1 from~\cite{he2024convexification}, we treat $\rho$ as a positive scaling variable. This allows us to derive two families of  valid linear inequalities for the nonlinear system in~\eqref{eq:transformation} with binary $\x$. The first  is obtained by scaling the standard McCormick relaxation for the bilinear terms $x_ix_j$:
% First, we use the variable $\rho$ to scale the following McCormick relaxation of $x_ix_j$~\citep{mccormick1976computability}:
\begin{equation*}
    0 \leq x_{i}x_j \leq x_i \quad \text{ and } \quad x_i+x_j-1 \leq x_{i}x_j  \leq x_j \qquad \for i \in [m] \text{ and } j \in [m]\setminus \{i\}. 
\end{equation*}
Using the definitions in~\eqref{eq:transformation}, this yields the following linear system in terms of $(\rho, \y, \z)$:
\begin{equation}\label{eq:mRMR_full_cons-4}
    0 \le z_{ij} \le y_i \quad \text{ and } \quad y_i + y_j - \rho \le z_{ij} \le y_j \qquad \for i \in [m] \text{ and } j \in [m]\setminus \{i\}.
\end{equation}
To derive the second class of inequalities, we consider a nonlinear representation of variable $\x$ given as follows:
\begin{equation}\label{eq:nl-x}
    x_k = \frac{x_k \cdot (\sum_{i,j \in [m]}x_ix_j)}{\sum_{i,j\in [m]}x_ix_j} \quad \for k \in [m],
\end{equation}
which holds since the denominator is strictly positive for all feasible $\x$. We next bound the numerator of~\eqref{eq:nl-x} from above and below by using the McCormick relaxation. For any $k \in [m]$, we have:
% Then, we use variable $\rho$ to scale the following over- and under-estimators of the numerator in~\eqref{eq:nl-x}:
\[
\begin{aligned}
    x_k \cdot \biggl(\sum_{i,j \in [m]}x_ix_j \biggr) 
&\leq   \min \biggl\{U^2x_k, L^2 x_k + \sum_{i,j \in [m]}x_ix_j - L^2 \biggr\}  , \\
x_k \cdot \biggl(\sum_{i,j \in [m]}x_ix_j \biggr) &\geq  \max \biggl\{L^2x_k, U^2x_k + \sum_{i,j \in [m]}x_ix_j - U^2  \biggr\},
\end{aligned}
\]
where $L^2$ (resp. $U^2$) is a valid lower (resp. upper) bound of the bilinear function $\sum_{i,j\in [m]}x_ix_j$ for $\x \in \X$. After scaling both nonlinear inequalities with $\rho$, and using the definitions in~\eqref{eq:transformation} and~\eqref{eq:nl-x}, we obtain the following system of linear inequalities in terms of variables $(\x, \rho,\y,\z)$: 
% \begin{equation}\label{eq:mRMR_full_cons-5}
%         \begin{aligned}
%             x_k &\leq \sum_{i\in [m]}z_{ik} +  \sum_{i\neq j}z_{ij} \quad \text{ and } \quad x_k \leq \sum_{i\in [m]}z_{ik} + (m^2-m)\cdot y_k \qquad  \for k \in [m]  \\
%     x_k &\geq \sum_{i\in [m]}z_{ik} \quad \text{ and } \quad x_k \geq \sum_{i\in [m]}z_{ik}  +  (m^2-m)\cdot (y_k-\rho) + \sum_{i\neq j}z_{ij} \qquad \for k \in [m].
%         \end{aligned}
% \end{equation}
\begin{equation}\label{eq:mRMR_full_cons-5}
        \begin{aligned}
            x_k &\leq U^2y_k \quad \text{ and } \quad x_k \leq L^2y_k + 1 - L^2\rho \qquad \for k \in [m] \\
    x_k &\geq  L^2y_k \quad \text{ and } \quad x_k \geq U^2y_k + 1 - U^2\rho \qquad \for k \in [m].
        \end{aligned}
\end{equation}
With inequalities in~\eqref{eq:mRMR_full_cons-4} and~\eqref{eq:mRMR_full_cons-5}, we are ready to present our MIP formulation for~\eqref{eqn: mRMR}: 
% \begin{equation}
%     L\rho \leq \sum_{i\in [m]} \leq U\rho \quad \text{ and } \quad 0\leq y_i \leq \rho \for i \in [m].
% \end{equation}
\begin{tightalign}
\begin{align}\label{eqn: mRMR_full_cons}   
        \max  \biggl\{ \sum_{i, j} \bigl(\MI(f_i, Y) - \MI(f_i, f_j)\bigr)  z_{ij} \biggm|
          \x \in \X,
        \sum_{i, j \in [m] } z_{ij} = 1,   z_{ii} = y_i \, \forall i \in [m], 
    \eqref{eq:mRMR_full_cons-4} \text{ and } \eqref{eq:mRMR_full_cons-5} \bigg\},\tag{\textsc{mRMR-Pers}}
\end{align}
\end{tightalign}
where $\X$ is the feasible region of~\eqref{eq:mrmr-fp}, the second constraint follows from the definition of $z_{ij}$ in~\eqref{eq:transformation}, and the third constraint is derived using the relation $x_ix_i = x_i$ for a binary variable. In Theorem~\ref{them:mRMR-Pers}, we show that~\eqref{eqn: mRMR_full_cons} is indeed a valid MIP formulation of~\eqref{eqn: mRMR}. 

Moreover, we present a theoretical comparison between our model and~\eqref{eq:mrmr-RMC}. One of the most important properties of an MIP formulation is the strength of its natural continuous relaxation that is obtained by ignoring the integrality constraints. This is important because a tighter continuous relaxation often indicates a faster convergence of the branch-and-bound algorithm on which most commercial MIP solvers are built~\citep{vielma2015mixed}. Since~\eqref{eqn: mRMR}  has a maximization objective, a tighter formulation has a smaller continuous relaxation objective value. Let $V_{\textsc{pers}}$ (resp. $V_{\textsc{rmc}}$) be the optimal objective value of the natural continuous relaxation of formulation~\eqref{eqn: mRMR_full_cons} (resp.~\eqref{eq:mrmr-RMC}).

\begin{theorem}\label{them:mRMR-Pers}
\eqref{eqn: mRMR_full_cons} is an MIP formulation of~\eqref{eq:mrmr-fp}. Moreover, $V_{\textsc{pers}} \leq V_{\textsc{rmc}}$. 
\end{theorem}

% Using the definitions in~\eqref{eq:transformation},~\eqref{eqn: mRMR} is equivalent to 
% \begin{equation}\label{eq:mrmr-trans}
%     \max  \biggl\{ \sum_{i, j} \bigl(\MI(f_i, Y) - \MI(f_i, f_j)\bigr)  z_{ij} \biggm|
%           \x \in \X,~\eqref{eq:transformation} \bigg\}.
% \end{equation}
% Since the objective functions of both problems,~\eqref{eq:mrmr-trans} and~\eqref{eqn: mRMR_full_cons}, are same, it suffices to show that their feasible regions are equivalent. By construction,~\eqref{eq:mRMR_full_cons-4} and~\eqref{eq:mRMR_full_cons-5} are valid inequalities for the feasible solution of~\eqref{eq:mrmr-trans}. Moreover, for every feasible solution $(\x,\rho,\y,\z)$ to~\eqref{eq:mrmr-trans}, it follows readily that $\sum_{i,j}z_{ij} =1$ and $z_{ii}= y_i$ for $i \in [m]$, where the second equality holds since $x_ix_i = x_i$. Thus, we obtain that the feasible region of~\eqref{eq:mrmr-trans} is contained in that of~\eqref{eqn: mRMR_full_cons}. To show the reverse, we consider a feasible solution $(\x,\rho,\y,\z)$ to~\eqref{eqn: mRMR_full_cons}. 

Last, we derive additional valid linear inequalities to tighten~\eqref{eqn: mRMR_full_cons}. Our inequalities are obtained by using the lower bound $L$ and upper bound $U$ on the total number of  selected features, and are given as follows: 
\begin{subequations}
 \label{eqn: mRMR_RLT}
 \begin{tightalign}
\begin{align}
         U y_i - \sum_{j \in [m]}z_{ij} &\geq 0   \quad \for i \in [m]\label{mRMR_RLT-1} \\
       \sum_{j \in [m]}z_{ij} - Ly_i &\geq 0 \quad \for i \in [m] \label{mRMR_RLT-2}
    \\ 
    U\rho - \sum_{j \in [m]} y_j - \biggl(U y_i - \sum_{j \in [m]}z_{ij} \biggr) &\geq 0
    % \quad v_{ii}^k = v_{ik}^k = v_{ki}^k = w_{ik}, 
   \quad  \for i \in [m]  \label{mRMR_RLT-3} \\
        \sum_{j \in [m]} y_j - L\rho - \biggl(\sum_{j \in [m]}z_{ij} - Ly_i\biggr) &\geq 0 \quad  \for i \in [m] \label{mRMR_RLT-4}  \\
         -\sum_{i,j \in [m]}z_{ij}+(U+L)\sum_{j \in [m]}y_j - (U\cdot L) \rho &\geq 0.  \label{mRMR_RLT-5}
 \end{align}
  \end{tightalign}
\end{subequations}
To derive these inequalities, we modify the reformulation-linearization technique (RLT) ~\citep{sherali1990hierarchy} as follows. First, we generate nonlinear inequalities by multiplying one of the bounding inequalities  with linear inequalities $0\leq x_i \leq 1$ for $i \in [m]$: 
\begin{equation*}
    \begin{aligned}
    \biggl(U - \sum_{j \in [m]}x_j\biggr) \cdot x_i \geq 0&  \quad \text{ and } \quad \biggl( \sum_{j \in [m]}x_j-L\biggr) \cdot x_i \geq 0,  \\
    \biggl(U - \sum_{j \in [m]}x_j\biggr) \cdot  (1- x_i ) \geq 0 & \quad \text{ and } \quad \biggl( \sum_{j \in [m]}x_j-L\biggr) \cdot (1- x_i) \geq 0. 
\end{aligned}
\end{equation*}
Then, scaling these nonlinear inequalities with $\rho$ yields inequalities~\eqref{mRMR_RLT-1}--\eqref{mRMR_RLT-4}. Similarly, the last inequality~\eqref{mRMR_RLT-5} is obtained by using $\rho$ to scale $(U-\sum_{i\in [m]}x_i)\cdot(\sum_{j\in [m]}x_j-L) \geq 0$. 

\subsection{Informative Review Subset Selection}\label{section:irssformulation}
Let $\x = (x_1, x_2, \ldots, x_n) \in \{0,1\}^n$ be a vector of binary decision variables modeling the subset of reviews selected, i.e., $x_i = 1$ if and only if the review $i$ is selected. Then, the main nonconvex component in~\eqref{eqn: Entropy_review_selection} can be expressed as 
\begin{tightalign}
\begin{align}\label{eq:logratio}
\log \Biggl(\frac{P_j\sum_{i \in [n]}x_i}{\sum_{i \in [n]}d_j^ix_i} \Biggr) \quad \for j \in [m],    \tag{\textsc{LogRatio}}
\end{align}
\end{tightalign}
where $\x \in \{0,1\}^n$ and $L \leq \sum_{i \in [n]}x_i \leq U$.  To treat this discrete nonconvex function, we  use the notion of convex extension. Given a function $f: S \subseteq \{0,1\}^n \to \mathbb{R}$, a continuous function $g(\cdot)$ is called convex (resp. concave) extension of $f(\cdot)$ if $g(\cdot)$ is convex (resp. concave) and $g(\x) = f(\x)$ for every $\x \in S$. Among infinitely many convex (resp. concave) extensions, the tightest one is called the convex (resp. concave) envelope of $f(\cdot)$, denoted as $\conv(f)(\cdot)$ (resp. $\conc(f)(\cdot)$). 

A key challenge in constructing extensions for~\eqref{eq:logratio} is that it is not defined at the origin $\boldsymbol{0} := (0,0,\dots,0)_{n\times 1}$.  To address this, we represent~\eqref{eq:logratio} as the difference of two composite functions defined on $\{0,1\}^n$, allowing us to construct extensions by treating each component separately. It is evident that our approach also yields an MIP formulation for the instance subset selection problem for the Naive Bayes classifier studied in~\cite{wei2015submodularity}, where the loss function also takes the form in~\eqref{eq:logratio}.

For each opinion $j \in [m]$, we introduce a pair of univariate outer-functions that extend the domain of the logarithm function to zero,
\begin{equation}\label{eq:univariate}
 \phi_j(z) := \begin{cases}
     \log(P_j \cdot z) &   z >0 \\
     2\log(P_j) - \log(2 P_j) & z = 0
 \end{cases} \quad \text{ and } \quad \psi_j(z) := 
 \begin{cases}
    \log(z) & z > 0 \\
            \phi_j(U) - \frac{\delta}{P_{j}}      & z = 0,
\end{cases} 
\end{equation}
where we recall that $P_j$ is the score of the $j^{\text{th}}$ opinion in the full data set,  $U$ is the upper bound on the total number of selected reviews, and $\delta$ is the penalty for not including opinion $j$, defined as in~\eqref{eq:penalty}. The values at 
$0$ are deliberately chosen  to ensure the validity of  Lemma~\ref{lem: review_selection_dod} and Propositions~\ref{prop: review_selection_submodular} and~\ref{prop:conc}, which are the building blocks of our final formulation~\eqref{eq:dataformulation}.  Based on~\eqref{eq:univariate}, we then define a pair of composite functions as follows:
\begin{equation}
f_j(\x) := \phi_j(\one^\intercal \x) \quad \text{ and } \quad g_j(\x) := \psi_j(\d^\intercal_j\x ) \quad \for \x \in \{0,1\}^n, 
\end{equation}
where $\one$ is the all-ones vector of dimension $n$ and for each opinion $j \in [m]$, $\d_j := (d^1_j,d^2_j, \ldots, d^n_j)$ indicates the presence of opinion $j$ in the dateset of reviews.  These definitions allow us to formulate~\eqref{eqn: Entropy_review_selection} as follows:
\begin{tightalign}
\begin{align}\label{eq:diff}
\begin{aligned}
    \min \biggl\{ \sum_{j \in [m]}P_j \bigl\vert f_j(\x)- g_j(\x) \bigr\vert \biggm| L \le\sum_{i \in [n]}x_i \le U,\  \x \in \{0, 1\}^n \biggr\}.
    % \text{s.t.} \quad  & \mu_j \geq \phi_j(\one^\intercal \x   ) - \psi_j(\d^\intercal_j\x ) \quad \text{ for } j \in [m] \\
    % & \mu_j \geq  \psi_j(\d^\intercal_j\x ) - \phi_j(\one^\intercal \x   ) \quad \text{ for } j \in [m]  \\
    % & L \le\sum_{i \in [n]}x_i \le U,\  
    % \x \in \{0, 1\}^n
\end{aligned} \tag{\textsc{RSKL-DC}}
 \end{align}
 \end{tightalign}
% More specifically, the function value of $\psi_j(\cdot)$ at $0$ allows us to model the penalty for not including opinion~$j$. This guarantees that~\eqref{eq:diff} is a reformulation of~\eqref{eqn: Entropy_review_selection}.  
\begin{lemma}
\label{lem: review_selection_dod}
Formulation~\eqref{eq:diff} is equivalent to formulation~\eqref{eqn: Entropy_review_selection}.
\end{lemma}
\begin{remark}
    Lemma~\ref{lem: review_selection_dod} allows us to focus on deriving polyhedral extensions for $f_j(\cdot)$  and $g_j(\cdot)$ separately, rather than working directly with the more complex function~\eqref{eq:logratio}. In Section~\ref{section:envelopes}, we exploit the specific structures of $f_j(\cdot)$ and $g_j(\cdot)$ to obtain their tight polyhedral extensions, which are then combined in Section~\ref{section:final} to yield an MIP formulation of~\eqref{eq:diff}.
\end{remark}

% Observe that, for each $j\in J$, the nonlinear expression inside the absolute value function can be expressed as differences of two functions, $\phi(P_{j}^D\sum_{i\in N} x_i ) - \psi_j(\sum_{i\in [n]}d_{ij}x_i)$, 

% The advantage of~\eqref{eq:diff} over~\eqref{eqn: Entropy_review_selection} is that *****   
\subsubsection{Explicit convex and concave envelopes.}\label{section:envelopes}
We begin by presenting the convex envelopes of $f_j(\cdot)$ and $g_j(\cdot)$. The function values at $\boldsymbol{0}$ are chosen specifically to ensure that  both functions are submodular. Recall that a function $f:\{0,1\}^n \to \R$ is submodular if
\begin{equation*}
   f(\x') + f(\x'')\geq f(\x' \vee \x'') + f(\x'\wedge \x'') \quad \for \x', \x'' \in \{0,1\}^n, 
\end{equation*}
where $\x' \vee \x''$ (resp. $\x' \wedge \x''$) is the componentwise maximum (resp. minimum) of $\x'$ and $\x''$. The submodularity of $f_j(\cdot)$ and $g_j(\cdot)$ allows us to describe  convex envelopes using the Lov\'asz extension of a set function~\citep{lovasz1983submodular}. The Lov\'asz extension is an extension of a function $f(\cdot)$ defined on $\{0,1\}^n$ to a function $f_{\Lovasz}(\cdot)$ defined on $[0,1]^n$. For each subset $S$ of $[n]$, let $\chi^S\in \{0,1\}^n$ be its indicator vector---that is, the $i^{\text{th}}$ coordinate of $\chi^S$ is $1$ if and only if $i \in S$. Observe that every vector $\x \in [0,1]^n$ can be expressed uniquely as $\x =  \lambda_0  \chi^{T_0} + \lambda_1  \chi^{T_1} + \cdots + \lambda_n \chi^{T_{n}} $,
where $\lambda_k \ge 0$ for $k = 0,1, \ldots,n $ and $\sum_{k =0}^n\lambda_k =1$, and $\emptyset = T_0 \subseteq T_1 \subseteq \cdots \subseteq T_n = [n]$. Thus, 
\begin{equation*}
    f_{\Lovasz}(\x) := \lambda_0 f(\chi^{T_0}) + \lambda_1 f(\chi^{T_1}) + \cdots + \lambda_n f(\chi^{T_n}) 
\end{equation*}
is a well-defined extension of the function $f(\cdot)$ (called the \textit{Lov\'asz extension} of $f(\cdot)$) on the continuous domain $[0,1]^n$. It is shown in~\cite{tawarmalani2013explicit} that the Lov\'asz extension $f_{\Lovasz}(\cdot)$ of $f(\cdot)$ coincides with its convex envelope if and only if $f(\cdot)$ is submodular. 
 Moreover, when the Lov\'asz extension $f_{\Lovasz}(\cdot)$ is convex, it is expressible as the pointwise maximum of affine functions~\citep{tawarmalani2013explicit}, that is,  $f_{\Lovasz}(\x) = \max_{\pi \in \Pi}\bigl\{f^\pi(\x)\bigr\}$ for every $\x \in [0,1]^n$,
% \[
% f_{\Lovasz}(\x) = \max_{\pi \in \Pi}\bigl\{f^\pi(\x)\bigr\} \quad \for \x \in [0,1]^n,
% \]
 where $\Pi$ is the set of all permutations of $[n]$ and $f^\pi(\cdot)$ is an affine function defined as follows:
% \[
%  f^\pi(\x):= \sum_{i=1}^{n} \biggl( f \Bigl(\sum_{j=1}^i e_{\pi(j)}\Bigr) - f\Bigl(\sum_{j=1}^{i-1}e_{\pi(j)}\Bigr)\biggr) \cdot x_{\pi(i)} +  f(\boldsymbol{0})  \quad \text{ for } \x \in [0,1]^n.
% \]
\[
 f^\pi(\x):= \sum_{i=1}^{n} \biggl( f \Bigl(\sum_{j=1}^i e_{\pi(j)}\Bigr) - f\Bigl(\sum_{j=1}^{i-1}e_{\pi(j)}\Bigr)\biggr) \cdot x_{\pi(i)} +  f(\boldsymbol{0})  \quad \text{ for } \x \in [0,1]^n.
\]
Now, we are ready to present explicit expressions of the convex envelope of $f_j(\cdot)$ and $g(\cdot)$.
% Under the definition of outer-functions $\phi_j(\cdot)$  and $\psi_j(\cdot)$ in~\eqref{eq:univariate}, the composite functions $f_j(\cdot)$ and $g_j(\cdot)$ are submodular and thus their convex envelopes can be obtained using the Lov\'asz extension. 
\begin{proposition}
\label{prop: review_selection_submodular}
Assume that $U \geq 2$. For each opinion $j \in [m]$,  $\conv(f_j)(\x) = (f_j)_{\Lovasz}(\x)$ and $\conv(g_j)(\x) =  (g_j)_{\Lovasz}(\x)$ for every $\x \in \{0,1\}^n$.
\end{proposition}

Next, we describe the concave envelope $f_j(\cdot)$ and $g_j(\cdot)$. Here, the function values at $\boldsymbol{0}$ enable us to treat both functions as compositions of univariate piecewise concave functions and linear functions with unit coefficients. This composition structure allows us to derive concave envelopes as follows. 
\begin{proposition}\label{prop:conc}
Assume that $U \geq 2$. For each $j \in [m]$, let $S_j:=\bigl\{i \in [n] \bigm| d_j^i = 1 \bigr\}$. Then
\begin{equation*}
\label{eqn: conc_f}
\conc(f_j)(\x) = \min_{k \in [n]}\biggl\{ \bigl[\phi_j(k  ) - \phi_j(k-1 )\bigr]\sum_{i\in [n]} x_i + k \phi_j\left(k-1\right)-(k-1)\phi_j\left(k\right) \biggr\} \quad \for \, \x \in [0,1]^n,
\end{equation*}
and 
\begin{equation*}
\label{eqn: conc_g}
    \conc(g_j)(\x) =\min_{k \in [\vert S_j \vert]} \biggl\{  \left[\psi_j(k) - \psi_j(k-1)\right]\sum_{i \in S_j} x_i + k \psi_j(k-1)-(k-1)\psi_j(k)  \biggr\} \quad \for \, \x \in [0,1]^n. 
\end{equation*}
\end{proposition}

\begin{remark}
Propositions~\ref{prop: review_selection_submodular} and~\ref{prop:conc} provide explicit convex and concave envelopes for the nonlinear functions $f_j(\cdot)$ and $g_j(\cdot)$. These envelopes allow us to replace $f_j(\cdot)$ and $g_j(\cdot)$ with their corresponding polyhedral extensions, which in turn lead to the final MIP formulation of~\eqref{eq:diff}.
\end{remark}

\subsubsection{Final formulation.}\label{section:final}
To obtain our formulation, we introduce additional variables to represent nonlinear structures. For each opinion $j \in [m]$, we introduce a variable $s_j$ (resp. $t_j$) to represent $f_j(\cdot)$ (resp. $g_j(\cdot)$), and a variable $\mu_j$ to present the absolute value function $\vert s_j - t_j\vert$. Using these additional variables, and replacing $f_j(\cdot)$ and $g_j(\cdot)$ with their envelope expressions, we arrive at an MIP formulation: 
\begin{tightalign}
\begin{align}
\label{eq:dataformulation}
    \begin{aligned}
            \min_{\x, \boldsymbol{\mu}, \boldsymbol{s},\boldsymbol{t}} \quad &  \sum_{j\in [m]}P_{j}\mu_j \\
    \text{s.t.} \quad &L \le\sum_{i \in [n]}x_i \le U  ,\  
    \x \in \{0, 1\}^n  \\
    &\mu_j \geq s_j - t_j  ,\
    \mu_j \geq -s_j + t_j  \quad \text{ for } j \in [m] \\
       &  s_j \geq ( f_j)_{\Lovasz}(\x) ,\  t_j \geq  ( g_j)_{\Lovasz}(\x) \quad \text{ for } j \in [m] \\
       & s_j \leq \conc(f_j)(\x) ,\ t_j \leq \conc(g_j)(\x) \quad \text{ for } j \in [m].
\end{aligned}\tag{\textsc{RSKL-Env}}
\end{align}
\end{tightalign}

\begin{theorem}\label{theorem:envelope}
Formulation~\eqref{eq:dataformulation} is an MIP formulation of Problem~\eqref{eqn: Entropy_review_selection}. 
\end{theorem}

\section{Implementations}\label{section:implementation}
% One of advantages of MIP formulations allows us to leverage the state-of-the-art solvers, such as \texttt{Gurobi}, 

% The MILP formulations in Section~\ref{sec: IP_formulation}, although tractable for small to medium instances, can be difficult to solve directly for large instances.  In branch-and-bound methods, a frequently observed source of inefficiency is that solvers explore highly suboptimal regions of the search space in considerable depth \citep{bertsimas2022scalable}. To mitigate this behavior, we adopt two strategies for the implementation of both formulations~(\ref{eqn: mRMR_full_cons}) and~\eqref{eq:dataformulation}. First, we supply a warm-start  (i.e., an initial feasible solution) that is installed as the incumbent by the solver, which could substantially tighten the initial lower bound of MILP formulations and guide the search toward regions with higher potential for optimality. Second, we incorporate optimality cuts to eliminate provably suboptimal regions.  These cuts allow the solver to prune inferior solutions early, thereby reducing the size of the feasible binary space explored in subsequent iterations and improving overall computational performance.

% Indeed, by pruning
% suboptimal solutions, it can encourage branch-and-bound methods to focus on regions of the search space that contain near-optimal solutions. 

In this section, we discuss the detailed framework for the implementations of formulations~\eqref{eqn: mRMR_full_cons} and~\eqref{eq:dataformulation} with warm-start and optimality cuts in Sections~\ref{sec: algorithm_mRMR} and~\ref{sec: algorithm_data}, respectively. Besides, for  formulation~\eqref{eq:dataformulation}, we integrate lazy fashion implementation with \texttt{Gurobi} \texttt{Callback} mechanism  to tackle the factorial number of constraints introduced by the Lov\'asz extension.

\subsection{Implementation of Formulation~\eqref{eqn: mRMR_full_cons}}
\label{sec: algorithm_mRMR}

% \cmt{Jun}{is there a better way to represent (1) and (2) in Algorithm 1  by other format?}

\SetInd{0.5em}{1em}
\SetNlSkip{1em}
Our implementation of formulation~\eqref{eqn: mRMR_full_cons} is shown in Algorithm~\ref{agm: mRMR_implementation}. In Step~1, we obtain a feasible solution via a backward elimination procedure proposed in~\cite{naghibi2014semidefinite}. Starting with a full set of features $S_1 = [m]$, backward elimination eliminates feature $\gamma \in S_j$ that maximizes $I_{\text{mRMR}}(S_j \setminus \{\gamma\})$ and updates $S_{j+1} = S_j \setminus \{\gamma\}$ at each iteration.  After iteration $m$, we select $S_{\text{init}}$ from $\{ S_j \mid L \le |S_j| \le U, j \in [m] \}$ that maximizes $I_{\text{mRMR}}(\cdot)$.  
%and use it as warm-start.
% Popular heuristics for mRMR optimization include forward selection, which starts with an empty set of features and adds one feature that contributes the most to mRMR performance, 
% and backward elimination, which starts with a full set of features and sequentially eliminates those that contribute the least to mRMR performance. 
% As backward elimination evaluate the contribution of a given feature in the context of all other features, it can consider the interaction among different features and show better performance compared with other heuristics according to \cite{naghibi2014semidefinite}. 
% Therefore, we assign the solution $\x_{\text{init}}$ obtained by backward elimination as the warm-start of variable $\x$ and use Equation~(\ref{eq:transformation}) to obtain the warm-start for variables $\rho$, $\y$, and $\z$. 
Second, we tighten formulation~\eqref{eqn: mRMR_full_cons} using an optimality cut defined as follows:
\begin{equation}
\label{eqn: mRMR_optimality_cut}
    V_{\text{init}} \le \sum_{i, j} \bigl(\MI(f_i, Y) - \MI(f_i, f_j)\bigr) \cdot z_{ij} \le V_{\text{relax}},
\end{equation}
where $V_{\text{init}} = I_{\text{mRMR}}(S_{\text{init}})$, 
% is the objective value of the feasible solution given by the backward elimination procedure, 
$V_{\text{relax}}$ is the optimal value of the LP relaxation of ~\eqref{eqn: mRMR_full_cons} tightened with the RLT constraints in~\eqref{eqn: mRMR_RLT}. After incorporating the RLT constraints and the optimality cut into~\eqref{eqn: mRMR_full_cons},  we call $\texttt{Gurobi}$ to solve the resulting formulation with $S_{\text{init}}$ being a warm-start solution. 

\begin{algorithm}[!t]
\linespread{1.2}\selectfont
\KwData{
Mutual information values $\MI(\gamma_i, Y)$ and $\MI(\gamma_i, \gamma_j)$ for $i,j \in [m]$, computed from the dataset%
\endnote{Continuous variables are discretized following \cite{brown2012conditional}, which is a standard treatment in information-theoretic feature selection algorithms \citep{li2017feature}. 
Mutual information is computed using the~\texttt{scikit-learn} package in~\texttt{Python} \citep{scikit-learn}.
}, and bounds $L$ and $U$ on the number of selected features.
}
\KwResult{An optimal solution to~(\ref{eqn: mRMR}).}

\textit{Step 1: Feasible solution  search  via backward elimination.}\\
\Indp Initialize $S_1 = [m]$ (full feature set)\;
At each iteration, remove the feature $\gamma \in S_j$ that maximizes $I_{\text{mRMR}}(S_j \setminus \{\gamma\})$, and update $S_{j+1} = S_j \setminus \{\gamma\}$, until $j = m$\;
From $\{ S_j \mid L \le |S_j| \le U, j \in [m] \}$, choose $S_{\text{init}}$ with the highest $I_{\text{mRMR}}(\cdot)$\;
\Indm

\textit{Step 2: Optimality cut generation.}\\
\Indp Set $V_{\text{init}} = I_{\text{mRMR}}(S_{\text{init}})$\;
Incorporate the RLT constraints~\eqref{eqn: mRMR_RLT} into~(\ref{eqn: mRMR_full_cons})\;
Solve the LP relaxation of~(\ref{eqn: mRMR_full_cons}) and record its optimal objective value as $V_{\text{relax}}$\;
Add the optimality cut~\eqref{eqn: mRMR_optimality_cut} to~(\ref{eqn: mRMR_full_cons})\;
\Indm

\textit{Step 3: Final optimization.}\\
\Indp Solve~(\ref{eqn: mRMR_full_cons}) using \texttt{Gurobi}, warm-starting from $S_{\text{init}}$\;
\Indm

\caption{Procedure for solving~(\ref{eqn: mRMR_full_cons}).}
\label{agm: mRMR_implementation}
\end{algorithm}
% by adding an optimality cut, which limits the lower bound and the upper bound of the optimal objective value. An immediate lower bound is the objective value of the feasible solution obtained by backward elimination, denoted as  $V_{\text{init}}$. Meanwhile, the upper bound, $V_{\text{relax}}$, corresponds to the optimal objective value of the continuous relaxation of formulation~\eqref{eqn: mRMR_full_cons}.  Then, the resulting optimality cut takes the following form:
% \begin{equation}
% \label{eqn: mRMR_optimality_cut}
%     V_{\text{init}} \le \sum_{i, j} \bigl(\MI(f_i, Y) - \MI(f_i, f_j)\bigr) \cdot z_{ij} \le V_{\text{relax}}.
% \end{equation}
% Since we have shown the tightness of the relaxation for formulation~(\ref{eqn: mRMR_full_cons}), (\ref{eqn: mRMR_optimality_cut}) can help accelerate the branch-and-bound method by enabling more effective pruning of the search space. 
% The detailed procedure is shown in Algorithm~\ref{agm: mRMR_implementation}. 

\subsection{Implementation of Formulation~\eqref{eq:dataformulation}}
\label{sec: algorithm_data}

\SetInd{0.5em}{1em}
\SetNlSkip{1em}
\SetAlgoNoLine
\begin{algorithm}[!t]
\linespread{1.2}\selectfont
\KwData{Review data matrix $D$,  bounds $L$ and $U$ on the number of selected instances.}
\KwResult{An optimal solution to~(\ref{eq:dataformulation}).}
\textit{Step 1: Feasible solution  search  via greedy  method.} \\ 
\Indp $S = \emptyset$\;
\For{$k\in [U]$}{Select $i\in [n]\setminus S$ that minimizes $I_{\text{KL}}(S\cup \{i\})$ and update $S = S\cup \{i\}$\;
}\
Set $S_{\text{init}} = S$\; \Indm
\textit{Step 2: Optimality cut generation.}\\
\Indp Set $V_{\text{init}} = I_{\text{KL}}(S_{\text{init}})$\;
Add~\eqref{eqn: review_opt_cut} as optimality cut to~\eqref{eq:dataformulation}\; \Indm
% Set $t = 1$ and $Ind = 1$\;
% \While{$Ind = 1$}{
% Solve $P_{t-1}$ and obtain current optimal solution $x^*$, $s^*$, $t^*$;\\
% Call Algorithm~\ref{agm: review_selection_seperation} to find the violated inequalities of current solution;\\
% \eIf{there is no violated inequalities}{Set $Ind = 0$}{Add the violated inequalities as cuts to $P_{t-1}$ and set the updated problem as $P_t$\;
% $t = t+1$\;}
\textit{Step 3: Final optimization.}\\
\Indp Solve~\eqref{eq:dataformulation} using the \texttt{Callback} routine of \texttt{Gurobi} with the separation oracle in Algorithm~\ref{agm: review_selection_seperation} and the warm-start solution $S_{\text{init}}$\; \Indm
\caption{Procedure for solving~\eqref{eq:dataformulation}}\label{agm: review_selection_cutting_plane}
\end{algorithm}

\SetInd{0.5em}{1em}
\SetNlSkip{1em}
\begin{algorithm}[!t]
\linespread{1.2}\selectfont
\KwData{Current solution $\x^*$, $\s^*$, $\t^*$, review data matrix $D$.} 
\KwResult{Violated inequalities.}
Generate $\pi \in \Pi([n])$ by sorting $\x^*$ such that $x^*_{\pi(1)} \geq x^*_{\pi(2)}\ge\cdots\ge x^*_{\pi(n)}$;\\
\For{$j\in [m]$}{
\If{$f_j(\x^*) - s^*_j   > \epsilon \cdot \vert f_j(\x^*) \vert$}{
% Generate $\pi \in \Pi([n])$ by sorting $\x^*$ such that $x^*_{\pi(1)} \geq x^*_{\pi(2)}\ge\cdots\ge x^*_{\pi(n)}$;\\
Return the violated inequality: $s_j \geq  f_j^\pi(\x)$;}

\If{$g_j(\x^*) - t^*_j  > \epsilon \cdot \vert g_j(\x^*) \vert$}{
% Generate $\pi \in \Pi(S_j)$ by sorting $\{x^*_i\}_{i\in S_j}$ such that $x^*_{\pi(1)} \ge x^*_{\pi(2)} \ge \cdots \ge x^*_{\pi(\vert S_j \vert)}$;\\
Return the violated inequality: $t_j \geq g_j^\pi(\x)$;
}
}
\caption{Separation oracle  for the  Lov\'{a}sz extension constraints}\label{agm: review_selection_seperation}
\end{algorithm}

In this subsection, we present our implementation of formulation~\eqref{eq:dataformulation} to solve large-scale Problem~\eqref{eqn: Entropy_review_selection}, as detailed in Algorithm~\ref{agm: review_selection_cutting_plane}. In Step~1, similar to Algorithm~\ref{agm: mRMR_implementation}, we provide an initial feasible solution with a greedy method, which is widely used in instance/review selection~\citep{lappas2012selecting, tsaparas2011selecting}. It starts with an empty set of reviews (i.e., $S = \emptyset$), and sequentially adds the review $i\in [n]\setminus S$ that minimizes $I_{\text{KL}}(S\cup \{i\})$ and updates $S = S \cup \{i\}$
% to the objective function in~\eqref{eqn: Entropy_review_selection} 
at each iteration until $\vert S\vert = U$. The resulting solution, denoted by $S_{\text{init}}$, can be used as a warm-start for \texttt{Gurobi}, and derive the following optimality cut 
\begin{equation}
\label{eqn: review_opt_cut}
     \sum_{j\in [m]}P_{j}\mu_j \le V_{\text{init}}, 
\end{equation}
where $V_{\text{init}}$ is the objective value of $S_{\text{init}}$. 

% One of the main difficulties of solving~\eqref{eq:dataformulation} is that the number of linear
% inequalities used to describe  Lov\'asz extensions is $n!$, which grows exponentially as the number of reviews $n$ becomes large. To address this, we implement~\eqref{eq:dataformulation}   in the so-called lazy fashion. More specifically, we use the \texttt{Callback} routine of ~\texttt{Gurobi} to implement Lov\'asz extensions as lazy constraints. Lov\'asz extensions are removed from~\eqref{eq:dataformulation} and placed in a lazy constraint pool. Then, \texttt{Gurobi} solves the relaxed problem and checks whether inequalities in the pool are violated at each integer solution generated in the branch-and-bound tree.  To help $\texttt{Gurobi}$ finding a violated inequality, we devise a separation procedure, Algorithm~\ref{agm: review_selection_seperation},  for the Lov\'asz extension constraints. For a given solution $(\x^*,\s^*,\t^*)$, we find a permutation $\pi$ such that $x^*_{\pi(1)} \geq \cdots \geq x^*_{\pi(n)}$. If  $f_j(\x^*) -s_j^* > \epsilon $ (resp.  $g_j(\x^*) -t_j^* > \epsilon $) for a small positive number $\epsilon$, we return $s_j \geq f^\pi_j(x)$ (resp. $t_j \geq g^\pi_j(x)$) as a violated inequality. 

One of the principal challenges in solving~\eqref{eq:dataformulation} is that the representation of Lov\'asz extensions requires $n!$ linear inequalities, whose number grows factorially with the size of the review set $n$. To mitigate this combinatorial explosion, we employ a lazy-constraint implementation. In particular, we utilize the \texttt{Callback} routine in \texttt{Gurobi} to enforce Lov\'asz extension inequalities dynamically. Rather than embedding these constraints directly in~\eqref{eq:dataformulation}, we exclude them from the initial formulation and instead maintain them in a lazy constraint pool. During the branch-and-bound procedure, \texttt{Gurobi} solves the relaxed problem and, at each integer feasible solution, verifies whether any inequalities from the pool are violated.

To facilitate this verification, we develop a fast separation routine (Algorithm~\ref{agm: review_selection_seperation}) for Lov\'asz extension inequalities. Given a candidate solution $(\x^*,\s^*,\t^*)$, we construct a permutation $\pi$ satisfying $x^*_{\pi(1)} \geq \cdots \geq x^*_{\pi(n)}$. If the condition $f_j(\x^*) - s_j^* > \epsilon \cdot \vert f_j(\x^*) \vert$ (resp. $g_j(\x^*) - t_j^* > \epsilon \cdot \vert g_j(\x^*) \vert$) holds for a prescribed tolerance $\epsilon > 0$, the routine identifies a violated inequality and returns the corresponding constraint $s_j \geq f_j^\pi(x)$ (resp. $t_j \geq g_j^\pi(x)$).

\section{Numerical Experiment}
\label{sec: numerical}

This section reports the computational performance of our MILP formulations for solving Problems~\eqref{eqn: mRMR} 
and~\eqref{eqn: Entropy_review_selection} in Sections~\ref{section:computation_mrmr} and~\ref{sec: result_RSKL}, respectively. All experiments were conducted on a personal laptop equipped with a 16-core, 32-thread AMD Ryzen 9 9950X CPU (base clock 4.6 GHz) and 96 GB of RAM. In all the experiments, we use
\texttt{Gurobi 12.0.1} as the optimization solver, within
the \texttt{Julia} programming language \citep{bezanson2017julia} with the \texttt{JuMP} modeling framework \citep{dunning2017jump}.  We set the time limit as 3600s and the optimality tolerance as 0.5\%.

\subsection{Feature Selection with mRMR}\label{section:computation_mrmr}
In this section, we evaluate four alternative formulations for solving Problem~\eqref{eqn: mRMR}, where the last three serve as benchmarks from the existing literature.

\begin{itemize}
    \item[-] \textsf{PersRLT}: our formulation~\eqref{eqn: mRMR_full_cons} implemented using Algorithm~\ref{agm: mRMR_implementation}.
    \item[-] \textsf{RMC}: formulation~\eqref{eq:mrmr-RMC}  with $\rho^U = \frac{1}{L^2}$ and $\rho^L = \frac{1}{U^2}$ in the implementation.
    \item[-] \textsf{BigM}: formulation ~\eqref{eq:mrmr-bigM} proposed by~\cite{nguyen2009optimizing}, which linearizes the trilinear terms with big-M technique as detailed in Appendix~\ref{apx: mRMR-bigm}. 
    \item[-] \textsf{VD}: formulation~\eqref{eq:mrmr-VD} proposed by~\cite{mehmanchi2021solving}, which treats the bilinear denominator with value-disjunction approach as detailed in Appendix~\ref{apx: mRMR-vd}. 
\end{itemize}
% \cmt{Jun}{shall we add one reference or a sentence somewhere to indicate the last three are benchmark algorithms?} 
We test the computational efficiency and the scalability of four formulations with both real and synthetic data. To ensure a fair comparison across alternative formulations in the numerical experiments, we initialize all formulations with the feasible solution obtained by the backward elimination.

\subsubsection{Experiment on real datasets.}
\label{sec: mRMR_real_data}

\begin{table}[!t]
  \centering
  \caption{Summary of real datasets.}
    \begin{tabular}{cccc}
    \hline
    Name  & \#Feature & \#Instance  & Source \\
    \hline
    \texttt{Statlog} & 19    & 2310 & \cite{singha2018adaptive}\\
        \texttt{Dermatology}& 34    & 366 & \cite{wan2022r2ci} \\
    \texttt{Lung Cancer} & 56    & 32 & \cite{naghibi2014semidefinite}\\
    \texttt{Optdigits} & 64 & 5620  & \cite{nguyen2014effective}\\ \texttt{Musk2}& 166& 6598 &\cite{gao2016variational}\\
    \texttt{Arrhythmia} & 278 & 370  & \cite{naghibi2014semidefinite}\\
    \texttt{Lung}    & 325   & 73 & \cite{gao2016variational}\\
    \texttt{GSE28700}& 556& 44 & \cite{wan2022r2ci} \\
    \texttt{Multiple Features} &649 & 2000 & \cite{nguyen2014effective}\\
    \texttt{CNAE-9 }    & 856   & 1080 & \cite{naghibi2014semidefinite}\\
    \hline
    \end{tabular}%
  \label{tab: dataset}%
\end{table}%
We begin by evaluating the computational efficiency of our formulation for~\eqref{eqn: mRMR} on real  datasets with varying dimensions and sizes, as summarized in Table~\ref{tab: dataset}. In this experiment, we fix $L = 1$ and $U=m$, that is, we consider all possible subsets of features. All datasets are obtained from the UCI Machine Learning Repository\endnote{See https://archive.ics.uci.edu.}, except for \texttt{GSE28700}, which is sourced from the Gene Expression Omnibus\endnote{See https://www.ncbi.nlm.nih.gov/geo/query/acc.cgi?acc=GSE28700.}. The corresponding URLs are provided in the endnotes for reproducibility and reference. The effectiveness of mRMR for feature selection on these datasets has been well-established in the literature, see details in Section~\ref{apx: mRMR_real_data}.

The performance metrics summarized in Table~\ref{tab: mRMR_real_data} provide a comprehensive assessment of computational efficiency and solution quality.
% \cmt{Jun}{is it possible to add one row about the performance of the proposed method in the reference paper to each dataset to make a fair comparison? Even in the original paper they may not report the time or root gap, at least the final gap is good enough.}
These include the total solution time for \texttt{Gurobi} in seconds (denoted as ``Time''), the number of nodes explored in the branch-and-bound search at termination (``Nodes''), the root gap, and the final optimality gap. The root gap, calculated as
\[
\text{Root gap} = \dfrac{|V_{\text{init}} - V_{\text{relax}}|}{|V_{\text{init}}|}\times 100\%,
\]
measures the relative difference between the objective value of the initial feasible solution ($V_{\text{init}}$) and the optimal objective value obtained from solving the root node relaxation ($V_{\text{relax}}$). A smaller root gap generally indicates a stronger initial solution or a tighter relaxation, both of which can significantly reduce overall computation time. The final gap, computed as 
\[
\text{Final gap} = \dfrac{|V_{\text{feas}} - V_{\text{bb}}|}{|V_{\text{feas}}|}\times 100\%,
\]
quantifies the relative difference between the best feasible solution identified by the solver ($V_{\text{feas}}$) and the best known upper bound ($V_{\text{bb}}$) at termination. A final gap of zero indicates that the solver has found a globally optimal solution, while a nonzero gap reflects the degree of remaining uncertainty in the solution quality.

\begin{table}[!t]
\TABLE
  {Comparison of~\textsf{PersRLT} with benchmark methods on real datasets.\label{tab: mRMR_real_data}}
  {  \begin{tabular}{cccrrrr}
    \hline
    Dataset & Features & \multicolumn{1}{c}{Method} & \multicolumn{1}{c}{Time (s)} & \multicolumn{1}{c}{Root gap (\%)} & \multicolumn{1}{c}{Final gap (\%)} & \multicolumn{1}{c}{Nodes} \\    \hline 
    \multirow{4}[0]{*}{ \texttt{Statlog}} & \multirow{4}[0]{*}{19} & 
  \textsf{\small PersRLT}    &  0.1     &     10.1  &    0.0  & 1 \\
         & & {\small \textsf{RMC} } &  2.2     &    73.6   &   0.0    &  4310 \\
        &  & \textsf{\small BigM} &  1.6   &   $\dagger$    &  0.0     & 203245 \\
        &  & \textsf{ \small VD} &   0.2    &     $\dagger$   &   0.0    &  2750\\ \hline 
          \multirow{4}[0]{*}{\texttt{Dermatology}} & \multirow{4}[0]{*}{34}& 
    \textsf{\small PersRLT}   &  0.4     &      16.1   &    0.0  & 242 \\
          & & {\small \textsf{RMC} } &  36.9     &   53.0  &   0.3    &  100216\\
          & & \textsf{\small BigM} &    3600   &   $\dagger$    &  71.1     & 38458259 \\
          & & \textsf{\small VD} &  1022.7    &     $\dagger$   &   0.0    &  763389\\ \hline 
    \multirow{4}[0]{*}{\texttt{Lung Cancer}} & \multirow{4}[0]{*}{56}&
    \textsf{\small PersRLT}   &  0.4     &      9.8  &    0.0  & 1 \\
          & & {\small \textsf{RMC} } &   2.7     &     152.4  &   0.0    &  2851\\
          & & \textsf{\small BigM} &    3600   &   $\dagger$    &  543.4     & 7362653 \\
          & & \textsf{\small VD} &   23.7    &     $\dagger$   &   0.0    &  880\\
          \hline
    \multirow{4}[0]{*}{\texttt{Optdigits}} &\multirow{4}[0]{*}{64}& 
    \textsf{\small PersRLT}   &   0.8   &   7.2    &    0.4   & 201 \\
          & & {\small \textsf{RMC} } &     3600   &   79.6    &  32.9     &  1619223\\
          & & \textsf{\small BigM} &    3600    &  $\dagger$     &  $\dagger$     &  36972840\\
          & & \textsf{\small VD} &    3600    &   $\dagger$    &    173.3   & 247885 \\
          \hline
    \multirow{4}[0]{*}{\texttt{Musk2}}&\multirow{4}[0]{*}{166} & 
    \textsf{\small PersRLT}   &    1.1   & 0.5      &  0.5     &  1\\
          & & {\small \textsf{RMC} } &     3600   &     74.2  &    51.6   &  294482\\
          & & \textsf{\small BigM} &    3600    &       102.4 &     95.9  & 7354131 \\
          & & \textsf{\small VD} &     3600   &   $\dagger$    &   75.7    &  2867\\
\hline
    \multirow{4}[0]{*}{\texttt{Arrhythmia}} &\multirow{4}[0]{*}{278}&
    \textsf{\small PersRLT}   &   84.1    &    17.3   &    0.0   & 173 \\
          & & {\small \textsf{RMC} } &    3600   &  162.4     &   101.5    &  931652\\
          & & \textsf{\small BigM} &     3600   &     $\dagger$  &    $\dagger$   &  6393247\\
          & & \textsf{\small VD} &       3600  &   $\dagger$    &   $\dagger$    & 13  \\
          \hline
    \multirow{4}[0]{*}{\texttt{Lung}} & \multirow{4}[0]{*}{325}&
   \textsf{\small PersRLT}   &  8.4     &    1.3   &     0.1  &  1\\
          & & {\small \textsf{\small RMC} } &      3600  &     11.4   &   6.0    & 49003\\
          & & \textsf{\small BigM} &     3600   &   $\dagger$    &   $\dagger$    & 5921779 \\
          & & \textsf{\small VD} &     3600    &    $\dagger$    &    $\dagger$   & 5826 \\
          \hline
    \multirow{4}[0]{*}{\texttt{GSE28700}} & \multirow{4}[0]{*}{556}&
   \textsf{\small PersRLT}   &   201.6    & 14.0      &   0.0    & 810 \\
          & & {\small \textsf{RMC} } &    3600   &   80.6   &   50.7    & 41199\\
          & & \textsf{\small BigM} &      3600  &   $\dagger$    &   $\dagger$    & 1509225 \\
          & & \textsf{\small VD} &      3600  &   $\dagger$    &  $\dagger$     &  1\\
          \hline
    \multirow{4}[0]{*}{\texttt{Multiple Features}} & \multirow{4}[0]{*}{649}&
   \textsf{\small PersRLT}  &     3600   &   10.6    &   6.0    &  5332\\
          & & {\small \textsf{RMC} } &   3600     &   69.8    &     54.3  & 3885\\
          & & \textsf{\small BigM} &     3600   &  $\dagger$     &   $\dagger$    & 1248056 \\
          & & \textsf{\small VD} &    3600    &     $\dagger$    &        $\dagger$ & 1 \\
          \hline
    \multirow{4}[0]{*}{\texttt{CNAE}} & \multirow{4}[0]{*}{856}&
   \textsf{\small PersRLT}   &   130.3    & 2.9      &    0.0   &  7 \\
          & & {\small \textsf{RMC} } &     3600   & 39.2      &   29.4    &  2206\\
          & & \textsf{\small BigM} &   3600    &  $\dagger$     &  $\dagger$     &  329145\\
          & & \textsf{\textsf{\small VD}} &   3600    &   $\dagger$     &  $\dagger$      & 1 \\
          \hline
    \end{tabular}}
{\textit{Note}. ``$\dagger$'' means the gap is lager than 1000\%. }
\end{table}%

The results in Table~\ref{tab: mRMR_real_data} show that our formulation~\textsf{PersRLT} successfully solves~\eqref{eqn: mRMR} to optimality on all datasets within the timelimit, except for~\texttt{Multiple Features}, where it terminates with a final optimality gap of 6\%. In contrast, the alternative formulations~\textsf{RMC},~\textsf{BigM}, and~\textsf{VD} fail to solve datasets that have more than 64 features, with final optimality gaps exceeding 50\% for~\texttt{Multiple Features}.
~\textsf{BigM}, and~\textsf{VD} have final optimality gaps that exceed 1, 000\% (as indicated by the "$\dagger$" symbol) for datasets with more than 200 features.

A key factor contributing to the poor performance of these benchmark methods is their weak root relaxation compared with~\textsf{PersRLT} as shown in Table~\ref{tab: mRMR_synthetic}. 
% The root gap of is substantially smaller than that of the three alternative formulations. 
In particular,~\textsf{PersRLT} consistently achieves a root gap below 20\% across all datasets, whereas~\textsf{BigM} and~\textsf{VD} exhibit root gaps exceeding 1,000\%. These results reveal the advantage of formulation~\textsf{PersRLT} in providing a significantly tighter root relaxation bound ($V_{\text{relax}}$) on the true integer optimal objective value, which improves computational efficiency by allowing for more effective pruning of the search tree and reducing the portion of the solution space that must be explored. As a result, our method~\textsf{PersRLT} converges to global optimality with less solution time and fewer branch-and-bound nodes compared to the benchmark methods.

\begin{table}[!t]
\TABLE
{Comparison of~\textsf{PersRLT} with benchmark methods on synthetic data.\label{tab: mRMR_synthetic}}
   { \begin{tabular}{cccrrrrrrrrr}
    \hline
\multirow{2}[0]{*}{m} & \multirow{2}[0]{*}{$[L, U]$} & \multirow{2}[0]{*}{Method}  & \multirow{2}[0]{*}{Sol} & \multicolumn{2}{c}{Time (s)}&\multicolumn{2}{c}{Root gap (\%)} & \multicolumn{2}{c}{Final gap (\%)} & \multicolumn{2}{c}{Nodes} \\
          &       &    &   &\multicolumn{1}{c}{Mean}  & \multicolumn{1}{c}{Std}   & \multicolumn{1}{c}{Mean}  & \multicolumn{1}{c}{Std}   & \multicolumn{1}{c}{Mean}  & \multicolumn{1}{c}{Std}   & \multicolumn{1}{c}{Mean} & \multicolumn{1}{c}{Std} \\
    \hline
    \multirow{8}[0]{*}{100} 
          & \multirow{4}[0]{*}{No} & \textsf{\small PersRLT} & 10 & 2.8  & 1.2       &43.5 & 19.1 & 0.0 & 0.0 &195   & 184\\
          &       & \small \textsf{RMC}    & 0  & 3600  & 0.6  & 161.0 & 76.8 & 33.3 & 21.5 & 325157  & 138695  \\
          &       & \textsf{\small BigM}  & 0  & 3600  & 0.2    & $\dagger$ & $\dagger$ & $\dagger$ & $\dagger$ & 7518542  & 900643  \\
          &       & \textsf{\small VD}  & 10   & 748.4  & 877.2      & $\dagger$ & $\dagger$ & 0.0 & 0.0 & 10757  & 12245  \\
              \cline{2-12}
          & \multirow{4}[0]{*}{\makecell{Yes}} & \textsf{\small PersRLT}   & 10     & 2.0  & 0.6    & 27.5 & 16.4 & 0.0 & 0.1 & 124  & 110    \\
          &       & \small \textsf{RMC}  & 1  & 3314.7  & 907.6     & 120.4 & 57.8 & 28.2 & 19.7 & 335813  & 140921  \\
          &       & \textsf{\small BigM}  & 0  & 3600  & 32.2    & $\dagger$ & $\dagger$ & $\dagger$ & 519.8 & 8809530  & 3585192  \\
          &       & \textsf{\small VD}  & 10   & 267.2  & 213.1      & $\dagger$ & $\dagger$ & 0.0 & 0.0 & 703   & 635  \\
          \hline
    \multirow{8}[0]{*}{300} 
          & \multirow{4}[0]{*}{No} &  \textsf{\small PersRLT}   & 10     & 486.9  & 519.6    & 16.9 & 3.1 & 0.0 & 0.0 & 10286  & 11325  \\
          &       & \small \textsf{RMC} & 0   & 3600  & 0.4     & 43.8 & 6.4 & 27.2 & 6.2 & 37943  & 4034  \\
          &       & \textsf{\small BigM}  & 0 & 3600  & 0.2     & $\dagger$ & $\dagger$ & $\dagger$ & $\dagger$ & 5912972  & 765506  \\
          &       & \textsf{\small VD}  & 0   & 3600  & 0.8     & $\dagger$ & $\dagger$ & $\dagger$& 136.4 & 1237  & 1034  \\
          \cline{2-12}

                    & \multirow{4}[0]{*}{\makecell{Yes}} & \textsf{\small PersRLT}   & 10     & 319.4  & 525.1     & 5.5 & 3.6 & 0.2 & 0.2 & 1100  & 2021  \\
          &       & \small \textsf{RMC}    & 0  & 3600  & 0.2   & 32.6 & 5.8 & 24.3 & 5.4 & 27420  & 7886  \\
          &       & \textsf{\small BigM}  & 0 & 3600  & 0.3     & $\dagger$ & $\dagger$ & $\dagger$ & 154.5 & 5175337  & 1012345  \\
          &       & \textsf{\small VD}  & 0  & 3600  & 0.4     & $\dagger$ & 712.8 & 270.6 & 76.6 & 35    & 33  \\

          \hline
    \multirow{8}[0]{*}{500} 
          & \multirow{4}[0]{*}{No} & \textsf{\small PersRLT}  & 7   & 1412.4  & 1538.1     & 12.6 & 2.6 & 2.1 & 3.2 & 5917  & 5578  \\
          &       & \small \textsf{RMC}  & 0  & 3600  & 0.8     & 32.9 & 3.6 & 24.3 & 4.5 & 4960  & 1370  \\
          &       & \textsf{\small BigM} & 0  & 3600  & 0.5     & $\dagger$ & $\dagger$ & $\dagger$ & $\dagger$ & 2766196  & 290128  \\
          &       & \textsf{\small VD}    & 0   & 3600  & 0.4   & $\dagger$ & $\dagger$ & $\dagger$ & $\dagger$ & 1  & 0  \\
          \cline{2-12}

                    & \multirow{4}[0]{*}{\makecell{Yes}} & \textsf{\small PersRLT} & 8     & 1279.9 & 1529.6      & 2.5 & 2.2 & 1.1 & 1.8 & 219   & 267  \\
          &       & \small \textsf{RMC}   & 0  & 3600  & 0.6    & 29.7 & 3.2 & 28.2 & 3.3 & 3298  & 822  \\
          &       &\textsf{\small BigM}  & 0  & 3600  & 0.8    & $\dagger$ & $\dagger$ & $\dagger$ & $\dagger$ & 3386480  & 292523  \\
          &       & \textsf{\small VD}   & 0   & 3600  & 0.5    & $\dagger$ & 491.7 & $\dagger$ & $\dagger$ & 1     & 0  \\

          \hline
    \multirow{8}[0]{*}{700} 
          & \multirow{4}[0]{*}{No} &  \textsf{\small PersRLT}  & 7   & 1591.9 & 1471.1      & 9.0 & 1.9 & 0.9 & 1.4 & 3782  & 4491  \\
          &       & \small \textsf{RMC}   & 0   & 3600  & 0.2    & 27.3 &4.4 & 20.5 & 3.9 & 1952  & 601  \\
          &       & \textsf{\small BigM}  & 0  & 3600  & 0.3    & $\dagger$ & $\dagger$ & $\dagger$ & $\dagger$ & 802494  & 257299  \\
          &       & \textsf{\small VD} & 0     & 3600  & 0.9  & $\dagger$ &$\dagger$ & $\dagger$ & $\dagger$ & 1  & 0  \\
          \cline{2-12}

                    & \multirow{4}[0]{*}{\makecell{Yes}} & \textsf{\small PersRLT}  & 8    & 1114.3  & 1551.8      & 1.3 & 1.6 & 1.1 & 1.5 & 57    & 150  \\
          &       & \small \textsf{RMC}   & 0  & 3600  & 0.2     & 34.1 & 3.2 & 33.7 & 3.2 & 756  & 80  \\
          &       &\textsf{\small BigM}  & 0 & 3600  & 1.0    & $\dagger$ & $\dagger$ & $\dagger$ & $\dagger$ & 687506  & 281298  \\
          &       & \textsf{\small VD}   & 0   & 3600  & 0.3    &$\dagger$ & 145.0 & $\dagger$& 290.1 & 1     & 0  \\
          \hline
    \end{tabular}}
{\textit{Note}. ``$\dagger$'' means the gap is lager than 1000\%. }
\end{table}

\subsubsection{Experiment on synthetic datasets.}
\label{sec: mRMR_synthetic}

The synthetic dataset with $n$ instances and $m$ features is generated following the procedure outlined in Section~8 of the online supplement from \cite{park2020subset} as detailed in Section~\ref{apx: mRMR_synthetic_data}. In this experiment, we fix $n = 30$ and generate 10 synthetic datasets for each $m\in \{100, 300, 500, 700\}$. To assess the impact of cardinality constraints, we consider two settings as indicated in column ``$[L, U]$'' of Table~\ref{tab: mRMR_synthetic}: (i) No cardinality constraint (denoted as “No”), where we set $L = 1$ and $U = m$; and (ii) With cardinality constraint (denoted as “Yes”), where we set $L = 0.1m$ and $U = 0.9m$. For each generated dataset, we solve the problem using~\textsf{PersRLT} as well as the benchmark methods. We report the mean and standard deviation of the total solution time, root gap, final gap, and the number of branch-and-bound nodes for each value of $m$ in Table~\ref{tab: mRMR_synthetic}. The column labeled ``Sol'' in Table~\ref{tab: mRMR_synthetic} indicates the number of datasets (out of 10) that were solved to optimality within the time limit.

% From Table~\ref{tab: mRMR_synthetic}, we can first see the time required to
% solve Formulation~\textsf{PersRLT} is consistently lowest of the four formulations across all scenarios with different number of features and cardinality constraints, which is aligned with experiment results for real data. For other three benchmark formulations, they cannot be solved to full optimality within the one-hour time limit when $m\geq 300$ for all of the instances while Formulation~\textsf{PersRLT} can be solved in most instances. Second, as the cardinality of the selected feature subset becomes less constrained, solving Formulation~\textsf{PersRLT} becomes more time-consuming and the root gap of  Formulation~\textsf{PersRLT} becomes wider.  This behavior is attributed to both the enlarged feasible region and the diminishing effectiveness of RLT-based relaxations under weaker cardinality constraints.

From Table~\ref{tab: mRMR_synthetic}, we observe that~\textsf{PersRLT} consistently achieves the lowest solution time among the four formulations across all scenarios, regardless of the number of features or the cardinality constraint settings. This observation is consistent with the experimental results on real-world datasets in Table~\ref{tab: mRMR_real_data}. In contrast, the other three benchmark formulations fail to solve any instances to optimality within the one-hour timelimit when $m\ge300$, whereas~\textsf{PersRLT} successfully solves the majority of instances in these settings. Additionally, as the cardinality constraint becomes less restrictive, solving~\textsf{PersRLT} requires more time and results in a wider root gap. This trend can be attributed to both the expansion of the feasible region and the reduced strength of the RLT-based relaxations when the cardinality bounds are loosened.

\subsection{Instance Selection with KL Divergence}
\label{sec: result_RSKL}
In this section, we first introduce the data generation method to generate synthetic review data. Then we evaluate the effectiveness and scalability of the proposed formulation~\textsf{Env} and compare it with baseline methods. To the best of our knowledge, no exact MIP formulation for problem~\eqref{eqn: Entropy_review_selection} has been proposed in the existing literature. Therefore, we use two heuristic methods from \citet{zhang2021review}—\textsf{ComS(1)} and \textsf{ComS(2)}—which are specifically designed to solve problem~\eqref{eqn: Entropy_review_selection}, as benchmarks.  The details of these methods are as follows.
\begin{itemize}
    \item[-] \textsf{Env}: our formulation~\eqref{eq:dataformulation} implemented using Algorithm~\ref{agm: review_selection_cutting_plane}, where we set $\epsilon = 0.01$.
    \item[-] \textsf{ComS($\theta$)}: This is the ComS heuristic method, which operates in two phases: It first relaxes the binary constraints and solve~\eqref{eqn: Entropy_review_selection} with nonlinear optimization solver. Then, it rounds to a feasible binary solution and expands the search space around the binary solution to obtain a better final solution. The search depth is controlled by the parameter $\theta$. In our experiment, we follow the  configurations in~\cite{zhang2021review} and set $\theta = 1$ (denoted as \textsf{ComS(1)})  and $\theta = 2$  (denoted as \textsf{ComS(2)}), respectively.
\end{itemize}

% First, we investigate the effect of data sparsity on~\textsf{Env}. Second, we compare~\textsf{Env} with existing heuristic approaches to demonstrate its superiority. Finally, we assess the scalability of~\textsf{Env} in solving large-scale review selection problems.

\subsubsection{Synthetic data generation.}
\begin{figure}[!t]
\centering
\begin{minipage}{0.45\linewidth}
        \centering
\includegraphics[width=0.92\columnwidth]{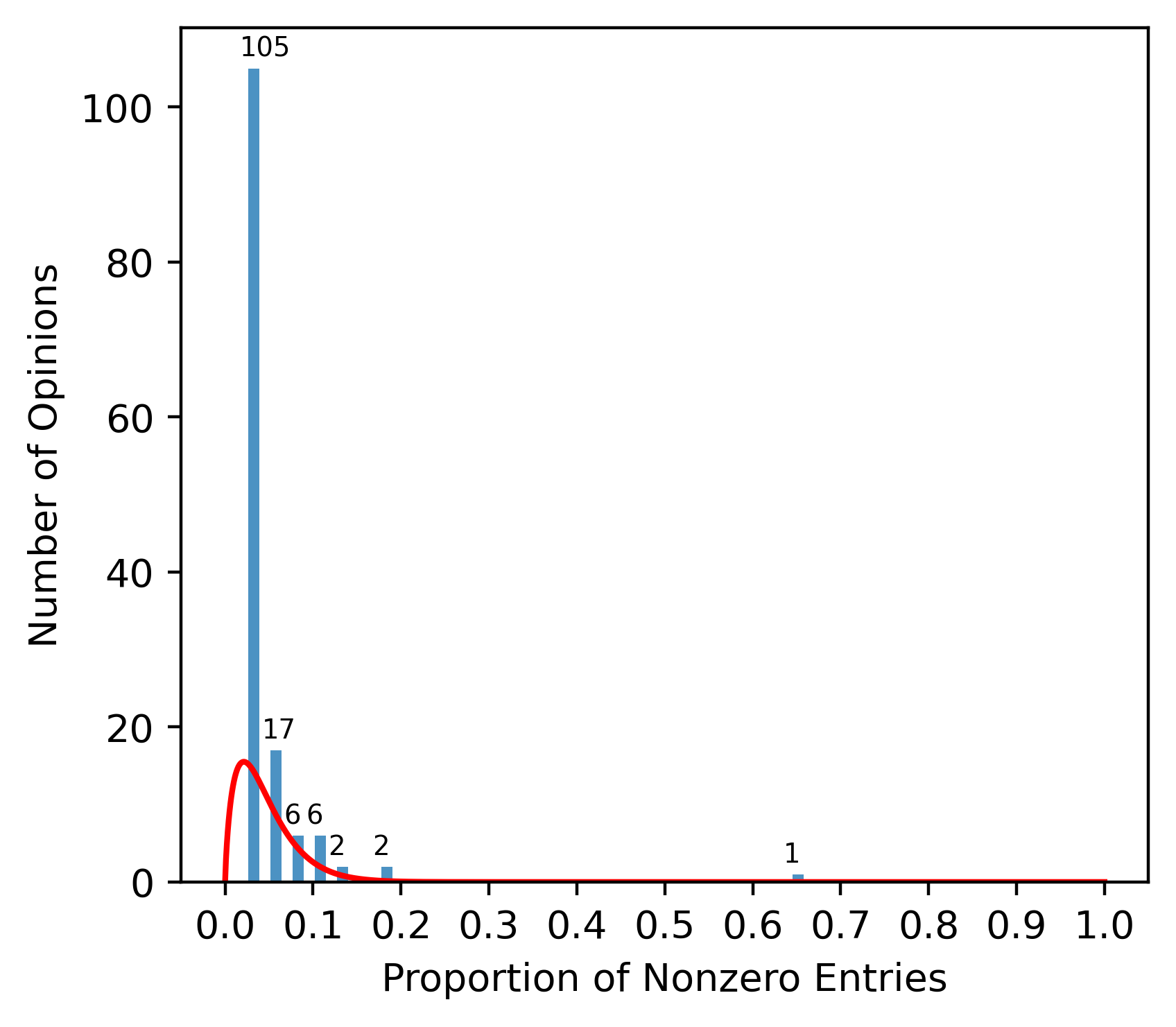}
		\caption{The histogram for the proportion of nonzeros entries of each opinion in sample reviews.}
		\label{fig: histogram}
\end{minipage}
\hspace{1cm}
\begin{minipage}{0.45\linewidth}
\centering
    \begin{tabular}{ >{\centering\arraybackslash}p{0.5cm}
    >{\raggedleft\arraybackslash}p{1cm}>{\raggedleft\arraybackslash}p{1cm}>{\raggedleft\arraybackslash}p{1.5cm}>{\raggedleft\arraybackslash}p{1cm}}
    \hline
    \multirow{1}[0]{*}{$\alpha$} & \multicolumn{1}{r}{\multirow{1}[0]{*}{$U$}}  & \multicolumn{1}{r}{\multirow{1}[0]{*}{Sol}} & \multicolumn{1}{r}{Time (s)} & \multicolumn{1}{c}{Nodes} \\
    \hline
    \multirow{3}[0]{*}{0.1} & 5      & 20   & 0.7      & 1   \\
          & 10  & 20   & 1.1     & 1     \\
          & 15   & 20   & 1.6     & 1     \\
          \hline
    \multirow{3}[0]{*}{1} & 5  & 20    & 0.7     & 1  \\
          & 10  & 20    & 1.1     & 1   \\
          & 15    & 20     & 1.7    & 1       \\
          \hline
    \multirow{3}[0]{*}{5} & 5   & 20     & 0.7    & 1      \\
          & 10     & 20   & 4.4    & 49   \\
          & 15   & 20   & 8.0     & 94      \\

          \hline
     \multirow{3}[0]{*}{10} & 5    & 20   & 1.4     & 17    \\
          & 10    & 20    & 25.0     & 211   \\
          & 15    & 20  & 505.0      & 686    \\
          \hline
    \end{tabular}%
\vspace{10pt}
\captionof{table}{Results of formulation~\textsf{ENV} for synthetic data with varied $\alpha$ ($\beta$ = 551, $n=200$, $m = 300$).}
\label{tab: review_selection_sparsity}
\end{minipage}  
\end{figure}

In the experiments of instance selection, the data are generated based on the sample reviews provided in Appendix~D of~\cite{zhang2021review}, which contains 38 reviews and 139 opinions for the camera sold on an online platform.
To analyze the distribution of opinion occurrences, we construct a histogram, shown in Figure~\ref{fig: histogram}.
% We first draw the histogram for the opinions of the sample reviews in Figure~\ref{fig: histogram} to analyze the distribution of the proportion of opinions occurrences. 
The horizontal axis denotes the proportion of nonzero entries in each opinion $j$, i.e., $P_j$, and the vertical axis indicates the number of opinions corresponding to a specific number of nonzero entry proportion. 
% \cmt{Jun}{it is better to describe the range of both horizontal and vertical axises. In addition, some of the vaules seem close to zero. will it be better to decribe it in table format instead of histogram?}
The histogram reveals a highly skewed distribution: the majority of opinions appear in a very small proportion of reviews. Specifically, 
% over 80\% of opinions occur in less than 5\% of the reviews, and 
more than 90\% of opinions appear in fewer than 10\% of the reviews. 
% It reveals that more than 90\% of opinions occur only once across the sample reviews, while more than 95\% of features exhibit fewer than 4 nonzero entries, i.e., less than 10\% of the total reviews 
% \textcolor{red}{need to rewrite}. 
This result underscores the high sparsity characteristic of real-world review data, reflecting the infrequent and uneven occurrence of opinions in user-generated content. 

% \begin{figure}[!t]
% \centering
% \includegraphics[width=0.8\textwidth]{histogram.png}
% \caption{The histogram for the proportion of nonzeros entries of each opinion ($P_j$) in sample reviews.}\label{fig: histogram}
% \end{figure}

Inspired by the sparse structure of the sample review data, we generate the synthetic data with $n$ observations and $m$ opinions in the following way: First, we calculate the proportion of nonzero entries for each opinion $j$ in the sample review data, denote it as $\hat P_{j}$. Second, we use beta distribution with shape parameters $\alpha$ and $\beta$ to fit the collection of proportions $\{\hat P_{j}\}_{j \in [139]}$ since $P_j$ is located in the interval $[0, 1]$,
and obtain the optimal shape parameters $\alpha = 0.12$ and $\beta = 551$\endnote{Beta distribution is fitted using~\texttt{SciPy} package in~\texttt{Python}.}. The fitted probability density function has been plotted in Figure~\ref{fig: histogram} with a red line. The low value of $\alpha$ combined with the high value of $\beta$ results in a beta distribution that is heavily skewed toward 0, capturing the sparsity pattern in the original data.
Third, we sample the proportions $\{P_j\}_{j\in [m]}$ from the fitted beta distribution, treating each $P_j$ as the probability that nonzero entries occur in opinion $j$. Finally, we generate each entry $d^i_j$ independently from a Bernoulli distribution with parameter~$P_j$, i.e., $d^i_{j}\sim Bern(P_j)$.

\subsubsection{Experiment on synthetic datasets.} 

In the following, we conduct three  experiments on the synthetic datasets to evaluate the performance of~\textsf{Env}.

In the first experiment, we investigate how the parameters $(\alpha, \beta)$ of the beta distribution affect the performance of~\textsf{Env}. We fix $n = 200$, $m = 300$, $L = 1$, and $\beta = 551$, while varying $\alpha \in \{0.1, 1, 5, 10\}$. For each value of $\alpha$, we randomly generate 20 synthetic instances. On each dataset, we apply~\textsf{Env} with different values of $U \in \{10, 20, 30\}$. The computational results are summarized in Table~\ref{tab: review_selection_sparsity}. We report the average final optimality gap, total running time, and number of branch-and-bound nodes. The column “Sol” indicates the number of instances (out of 20) solved within the time limit. 

From Table~\ref{tab: review_selection_sparsity}, we draw two main observations.
First, when the data is sparse (e.g., $\alpha \in \{0.1, 1\}$), the running times of~\textsf{Env} are negligible. However, as the data becomes denser (e.g., $\alpha \in \{5, 10\}$), the running times increase substantially. This trend arises because denser data contains more nonzero entries $d_j^i$, leading to a greater number of potential review subsets with nearly optimal objective values. As a result, the search for the optimal subset becomes more computationally demanding. To strike a balance between avoiding trivial sparsity and capturing realistic scenarios, we use $\alpha = 5$ and $\beta = 551$ for generating review data in the subsequent experiments. Second, as the upper bound $U$ on the cardinality of the selected review subset increases from 10 to 30, the problem becomes more computationally challenging.

In the second experiment, we compare~\textsf{Env} with~\textsf{ComS(1)} and~\textsf{ComS(2)}. We fix $m = 50$, $L = 1$, and vary $n \in \{50, 100, 150, 200\}$ following~\cite{zhang2012unified}. For each $n$, we randomly generate 20 review datasets and solve~\eqref{eqn: Entropy_review_selection} using the three methods with $U \in \{5, 10, 15\}$. The results are reported in Table~\ref{tab: review_selection_comparison}. For~\textsf{Env}, we present the mean and standard deviation of the final gap and running time for each $(n, U)$ combination. For~\textsf{ComS(1)} and~\textsf{ComS(2)}, we report the mean and standard deviation of running times and the relative gap (``Rel gap''), defined as
\begin{equation*}
    \text{Rel gap} = \frac{D_{\text{KL}}^{\textsf{ComS}}-D_{\text{KL}}^*}{D^{*}_{\text{KL}}} \times 100\%,
\end{equation*}
where $D_{\text{KL}}^*$ denotes the objective value of the solution obtained by~\textsf{Env}, and $D_{\text{KL}}^{\textsf{ComS}}$ is the objective value from either \textsf{ComS(1)} or \textsf{ComS(2)}.

\begin{table}[!t]

  \centering
  \caption{Comparison for methods~\textsf{ENV} and~\textsf{ComS($\theta$)}.}
    \begin{tabular}{crrrrrrrrrrrrr}
    \hline
    \multirow{3}[0]{*}{$n$} & \multicolumn{1}{c}{\multirow{3}[0]{*}{$U$}} & \multicolumn{4}{c}{\textsf{ENV}} & \multicolumn{4}{c}{\textsf{ComS(1)}}   & \multicolumn{4}{c}{\textsf{ComS(2)}} \\
    \cline{3-14}
          &    &   \multicolumn{2}{c}{Final gap (\%)} & \multicolumn{2}{c}{Time (s)} & \multicolumn{2}{c}{Rel gap (\%)} & \multicolumn{2}{c}{Time (s)} & \multicolumn{2}{c}{Rel gap (\%)} & \multicolumn{2}{c}{Time (s)} \\
          
          &    &   \multicolumn{1}{c}{Mean} & \multicolumn{1}{c}{Std} & \multicolumn{1}{c}{Mean} & \multicolumn{1}{c}{Std} & \multicolumn{1}{c}{Mean} & \multicolumn{1}{c}{Std} & \multicolumn{1}{c}{Mean} & \multicolumn{1}{c}{Std} & \multicolumn{1}{c}{Mean} & \multicolumn{1}{c}{Std} & \multicolumn{1}{c}{Mean} & \multicolumn{1}{c}{Std} \\
          \hline
    \multirow{3}[0]{*}{50} & 5     & 0.1 & 0.1 & 0.1  & 0.5  & 28.2  & 8.1  & 0.1  & 0.2  & 25.3  & 7.6  & 1.8  & 0.8  \\
          & 10  &  0.2 & 0.1  & 0.1  & 0.0  & 54.1  & 14.3  & 0.4  & 0.4  & 49.3  & 13.9 & 5.5  & 1.5  \\
          & 15  &0.2 & 0.1  & 0.1  & 0.0  & 99.2  & 33.0  & 0.5  & 0.4  & 91.2  & 31.6  & 11.1  & 2.9  \\
          \hline
    \multirow{3}[0]{*}{100} & 5     & 0.1 & 0.0  & 0.1  & 0.1  & 20.0  & 4.6  & 0.5  & 0.5  & 19.2 & 4.8  & 10.6  & 2.6  \\
          & 10 &  0.0 & 0.0  & 0.1  & 0.0  & 43.5  & 9.5  & 1.8  & 1.8  &41.5  & 7.9  & 44.0  & 11.0  \\
          & 15  & 0.0 & 0.0  & 0.1  & 0.1  & 70.8  & 18.2  & 6.3  & 3.9  & 65.4  & 19.9  & 99.0  & 30.2  \\
          \hline
    \multirow{3}[0]{*}{150} & 5 &  0.0  & 0.0  & 0.1  & 0.1 & 17.6  & 5.3  & 1.3  & 1.3  & 17.2  & 5.2  & 35.1  & 11.0  \\
          & 10    & 0.1 & 0.0  & 5.1  & 22.6  & 41.8  & 14.1  & 3.3  & 2.6  & 39.6  & 13.1  & 164.3  & 73.6  \\
          & 15   & 0.0 & 0.0  & 3.0  & 7.4  & 59.4  & 19.4  & 10.5  & 12.8  & 53.7  & 16.6  & 337.7  & 64.3  \\
          \hline
    \multirow{3}[0]{*}{200} & 5   & 0.0  & 0.0  & 0.2  & 0.3  & 20.2  & 4.5  & 0.7  & 0.7  & 19.3  & 4.0  & 267.6  & 363.5  \\
          & 10   & 0.0 & 0.0  & 2.7  & 5.8  & 47.5  & 9.5  & 2.7  & 3.9  & 44.2  & 7.3  & 707.0  & 752.1  \\
          & 15  & 0.0 & 0.0 & 25.4  & 81.5  & 77.1  & 18.9  & 2.7  & 5.1  & 72.0  & 19.2  & 1048.6  & 842.8  \\
          \hline
    \end{tabular}%
  \label{tab: review_selection_comparison}
\end{table}%

From Table~\ref{tab: review_selection_comparison}, two key findings emerge.
First, \textsf{Env} consistently achieves superior solution quality: the relative gaps of both \textsf{ComS} methods remain above 17\% across all settings, and the gaps widen as $U$ increases, indicating that~\textsf{ComS} performs worse for larger instance selection problems.
Second, the computation times of~\textsf{Env} are comparable to~\textsf{ComS(1)} and substantially lower than~\textsf{ComS(2)}—in fact, more than 40 times faster when $n = 200$. This efficiency gap reflects the theoretical time complexities: $O(mn^2U)$ for~\textsf{ComS(1)} and $O(mn^3U^2)$ for~\textsf{ComS(2)} \citep{zhang2021review}, which makes~\textsf{ComS(2)} computationally prohibitive for larger instances. Overall, these results highlight the effectiveness of formulation~\textsf{Env} in producing high-quality solutions to~\eqref{eqn: Entropy_review_selection} within practical running times.

In the final experiment, we examine the scalability of~\textsf{Env} in solving large-scale instance selection problems. We set $m = 300$, $L = 1$, and vary the number of reviews $n \in \{100, 200, \ldots, 700\}$. For each $n$, we randomly generate 20 review datasets and solve~\eqref{eqn: Entropy_review_selection} using~\textsf{Env} with $U$ set to $5\%$, $10\%$, and $15\%$ of $n$, so that the subset size scales proportionally with $n$. The results, summarized in Table~\ref{tab: review_selection_largescale}, include the number of solved instances (out of 20), the mean, standard deviation, and maximum final gaps, as well as the mean and standard deviation of running times, branch-and-bound nodes, and the number of constraints added during the callback procedure (denoted as “Lazy constraints”).

Table~\ref{tab: review_selection_largescale} provides two main insights.
First, as $n$ increases, both solution times and branch-and-bound nodes grow, reflecting the increased difficulty of larger instances. This trend stems from the higher dimensionality of the integer decision space and the expanding number of constraints, particularly those arising from the convex and concave envelopes in~\eqref{eq:dataformulation}. Therefore, while~\textsf{Env} is able to handle a wide range of large-scale settings, it exhausts the one-hour time budget before certifying optimality in the most challenging case ($n=700$, $U=0.15n$).
Second, despite the factorial growth in the number of potential constraints from the Lov\'{a}sz extension,~\textsf{Env} effectively manages this challenge. By employing lazy constraint generation, the algorithm activates only the necessary constraints, enabling the efficient solution of high-dimensional instances.

\begin{table}[!t]
  \centering
  \caption{Performance of Formulation~\textsf{Env} for large-scale review data ($\alpha = 5$, $\beta = 551$).}
    \begin{tabular}{crrrrrrrrrrr}
    \hline
    \multirow{2}[0]{*}{$n$ ($m = 300$)} & \multicolumn{1}{c}{\multirow{2}[0]{*}{$U$ (\%)}} & \multicolumn{1}{c}{\multirow{2}[0]{*}{Sol}} & \multicolumn{3}{c}{Final gap (\%)} & \multicolumn{2}{c}{Time (s)}  & \multicolumn{2}{c}{Nodes} & \multicolumn{2}{c}{Lazy constraints} 
    \\ 
          &       & & \multicolumn{1}{c}{Mean} & \multicolumn{1}{c}{Std} & \multicolumn{1}{c}{Max} & \multicolumn{1}{c}{Mean} & \multicolumn{1}{c}{Std} &   \multicolumn{1}{c}{Mean} & \multicolumn{1}{c}{Std} & \multicolumn{1}{c}{Mean} & \multicolumn{1}{c}{Std} \\ 
          \hline
    \multirow{3}[0]{*}{100} & 5   & 20   & 0.4  & 0.1  & 0.4  & 0.3  & 0.5   & 1     & 0     & 892  & 9  \\
          & 10  & 20  & 0.3   & 0.2  & 0.5  & 0.3  & 0.0      & 1     & 0     & 948   & 213  \\
          & 15  & 20  & 0.1  & 0.2  & 0.5  & 0.5  & 0.2    & 6     & 15    & 1398  & 1506  \\
          \hline
    \multirow{3}[0]{*}{200} & 5    & 20   & 0.2  & 0.2  & 0.4  & 0.8  & 0.1     & 1     & 0     & 905   & 14  \\
          & 10   & 20   & 0.0  & 0.1  & 0.3  & 7.6  & 11.1    & 51    & 71    & 8730  & 10188  \\
          & 15    & 20    & 0.0  & 0.1  & 0.3  & 16.1  & 14.3   & 110   & 89    & 16753  & 11010  \\
          \hline
    \multirow{3}[0]{*}{300} & 5   & 20    & 0.3  & 0.2  & 0.5  & 1.6  & 0.5     & 3     & 8    & 1323  & 1753  \\
          & 10 & 20    & 0.1  & 0.1  & 0.4  & 24.5  & 21.7      & 165   & 158   & 15469  & 9591  \\
          & 15   & 20   & 0.1  & 0.1  & 0.5  & 69.5  & 72.8    & 183   & 259   & 22863  & 15960  \\
          \hline
    \multirow{3}[0]{*}{400} & 5   & 20     & 0.1  & 0.2  & 0.5  & 7.4  & 9.4    & 58    & 142   & 7364  & 9054  \\
          & 10   & 20 & 0.2  & 0.2  & 0.5  & 193.4  & 411.5       & 376   & 544   & 37428  & 35964  \\
          & 15   & 15  & 0.5  & 0.9  & 3.8  & 1611.5  & 1632.2      & 865   & 679   & 39968  & 45530  \\
          \hline
    \multirow{3}[0]{*}{500} & 5   & 20    & 0.1  & 0.2  & 0.4  & 59.6  & 122.0     & 100   & 169   & 13182  & 7578  \\
          & 10  & 18    & 0.4  & 0.3  & 1.1  & 792.0  & 1083.8     & 706   & 539   & 77059  & 65405  \\
          & 15   & 16    & 0.6  & 1.2  & 5.0  & 1106.5  & 1350.8    & 795   & 622   & 64328  & 41894  \\
          \hline
    \multirow{3}[0]{*}{600} & 5   & 20    & 0.2  & 0.2  & 0.5  & 75.9  & 98.5     & 174   & 195   & 24311  & 20315  \\
          & 10   & 7    & 2.2  & 2.3  & 8.3  & 2922.1  & 1104.9   & 984   & 347   & 142780  & 40054  \\
          & 15   & 7    & 4.6  & 5.3  & 17.9  & 2643.7  & 1409.0   & 1168  & 571   & 105803  & 56424  \\

          \hline
    \multirow{3}[0]{*}{700} & 5  & 20      & 0.3  & 0.2  & 0.5  & 482.5  & 642.6    & 407   & 286   & 68641  & 49645  \\
          & 10  & 5   & 4.0  & 3.9  & 11.3  & 3068.6  & 1157.2     & 963   & 306   & 171172  & 63871  \\
          & 15   & 0  & 95.9  & 9.0  & 100.0  & 3600  & 5.7     & 605   & 569   & 65103  & 41153  \\

          \hline
    \end{tabular}%
  \label{tab: review_selection_largescale}%
\end{table}%

\section{Conclusion}\label{section:conlcusion}

% In this paper, we address the challenge of optimally solving information-theoretic data reduction by proposing novel MIP formulations for two representative problems: feature selection under the mRMR criterion and \replace{review selection}{instance selection} with KL divergence. For both problems, we develop exact and compact formulations derived from new polyhedral relaxations that capture the underlying nonlinear structures—bilinear and fractional terms in mRMR, as well as log-rational terms in KL divergence, which represent a broad class of information-theoretic objectives commonly encountered in data reduction. Extensive computational experiments on both synthetic and real-world datasets demonstrate that our formulations consistently help produce high-quality solutions for practically sized instances within reasonable time limits, significantly outperforming existing benchmark methods. Future study \textcolor{red}{xxx}

In this paper, we tackle the challenge of optimally solving two information-theoretic data reduction problems: feature selection under the mRMR criterion and instance selection using KL divergence. For both cases, we derive exact MIP formulations based on  polyhedral relaxations for nonlinear structures—bilinear and fractional terms in mRMR, and log-rational terms in KL divergence.  Through extensive experiments on both synthetic and real-world datasets, we show that our approach consistently yields high-quality solutions for practically sized instances within a reasonable time, substantially outperforming existing benchmark methods.

%\THEEndNotes
\begingroup \parindent 0pt \parskip 0.0ex \def\enotesize{\normalsize} \theendnotes \endgroup

% Appendix here
% Options are (1) APPENDIX (with or without general title) or
%             (2) APPENDICES (if it has more than one unrelated sections)
% Outcomment the appropriate case if necessary
%
% \begin{APPENDIX}{<Title of the Appendix>}
% \end{APPENDIX}
%
%   or
%
% \cmt{Jun}{you may add \{\} to capitalize the word in the reference title, e.g., {PCA} in Dey et al. (2022).}

% \bibliographystyle{apalike}
\bibliographystyle{informs2014}
\bibliography{reference}

\clearpage

\ECSwitch

\ECDisclaimer
\ECHead{Appendix}

\section{Supplementary Experimental Details for Section~\ref{section:computation_mrmr}}

\vspace{0.5cm}

\subsection{Alternative Benchmark MIP Formulations of~\eqref{eqn: mRMR}}\label{app:MIPmRMR}

In this section, we present the  benchmark MIP formulations used in Section~\ref{section:computation_mrmr} for globally solving Problem~\eqref{eqn: mRMR} with big-M and value-disjunction approaches. For more details, interested readers may refer to~\cite{nguyen2009optimizing} and~\cite{mehmanchi2021solving}.

\subsubsection{Big-M formulation.}
\label{apx: mRMR-bigm}
Let $c_{jk} = \MI(\gamma_j, Y) - \MI(\gamma_j, \gamma_k)$ for $j,k\in[m]$, then~\eqref{eq:mrmr-fp} can be rewritten as 
\begin{tightalign}
    \begin{align}
    \begin{aligned}
        \max_{} \quad &\sum_{j\in [m]} \biggl ( \sum_{k\in [m]} c_{jk}x_k \rho\biggr)x_j\\
         \text{s.t.} \quad &\sum_{j\in [m]} \biggl ( \sum_{k\in [m]} x_k \rho\biggr)x_j = 1, \\
        &  L \le \sum_{k\in [m]}  x_k \le U, \quad \x \in \{0, 1\}^m.
         \end{aligned}
    \end{align}
\end{tightalign}
Let $v_j^b := \biggl ( \sum_{k\in [m]} c_{jk}x_k \rho\biggr)x_j$ and $v_j^d := \biggl ( \sum_{k\in [m]} x_k \rho\biggr)x_j$,~\cite{nguyen2009optimizing} first use the big-M technique to linearize the product of $\sum_{k\in [m]} c_{jk}x_k \rho$ (resp. $\sum_{k\in [m]} x_k \rho$) and $x_j$ for $v_j^b$ (resp. $v_j^d$), and then use the McCormick envelope to linearize the bilinear term $y_k := x_k \rho$ in $\sum_{k\in [m]} c_{jk}x_k \rho$ and $\sum_{k\in [m]} x_k \rho$. The final formulation of~\eqref{eqn: mRMR} becomes
\begin{tightalign}
\begin{align}\label{eq:mrmr-bigM}
\begin{aligned}
  \max_{} \quad &\sum_{j\in [m]}v_j^b\\
   \text{s.t.} \quad &\sum_{j\in [m]}v_j^d = 1\\
  & -M_j^bx_j \leq v_j^b \leq  M_j^b x_j && \for j\in [m]\\
   &M_j^b(x_j - 1) + \sum_{k\in [m]}c_{jk}y_k \leq v_j^b \leq  M_j^b (1-x_j)+\sum_{k\in [m]}c_{jk}y_k && \for j\in [m]\\
   & -M_j^dx_j \leq v_j^d \leq  M_j^d x_j && \for j\in [m]\\
   & M_j^d(x_j - 1) + \sum_{k\in [m]}y_k \leq v_j^d \leq  M_j^d (1-x_j)+\sum_{k\in [m]}y_k && \for j\in [m]\\
   & 0\leq y_j\leq \rho^U x_j, \quad \rho^U(x_j - 1) + \rho \leq y_j \leq \rho && \for j\in [m]\\
   &  L \le \sum_{k\in [m]}  x_k \le U, \quad \x \in \{0, 1\}^m.
\end{aligned}
% \tag{\textsc{mRMR-BigM}}
\end{align} 
\end{tightalign}
% where $v_j^b$ replaces $\sum_{k \in [m]} c_{jk} z_{jk}$ and $v_j^d$  replaces $\sum_{k \in [m]} z_{jk}$ in~\eqref{mrmr-tri}.
We set $M_j^b = \sum_{k\in [m]}\vert c_{jk}\vert$, $M_j^d = m$, $\rho^U = 1$ following \cite{mehmanchi2021solving}.

\subsubsection{Value-disjunction (VD) formulation.}
\label{apx: mRMR-vd}
Let $c_{jk} = \MI(\gamma_j, Y) - \MI(\gamma_j, \gamma_k)$ for $j,k\in[m]$. Observing that $\sum_{j\in [m]}\sum_{k\in [m]}x_j x_k$ takes values in $\{1^2, 2^2, \ldots, m^2\}$,~\cite{mehmanchi2021solving} first reformulate~\eqref{eq:mrmr-fp} with the value-disjunction approach as follows:
\begin{tightalign}
\begin{align}
\begin{aligned}
  \max \quad &\sum_{l\in [m]}\sum_{j\in [m]}\sum_{k\in [m]}\frac{c_{jk}w_l x_j x_k}{l^2}\\
   \text{s.t.} \quad  & \sum_{j\in [m]}x_j = \sum_{l\in [m]}lw_l, \quad \sum_{l\in [m]}w_l = 1, \\
   &\quad L \le \sum_{k \in [m]} x_k \le U, \quad \x, \boldsymbol{w} \in \{0, 1\}^m.
   \end{aligned}
% \tag{\textsc{mRMR-VD}}
\end{align} 
\end{tightalign}
Then let $r:= \sum_{j\in [m]}\sum_{k\in [m]} c_{jk}x_j x_k$ and $s_l := r w_l$, ~\cite{mehmanchi2021solving} use big-M technique to linearize $r w_l$ and use McCormick envelope to linearize $t_{jk} := x_j x_k$, which results in the following final formulation:
\begin{tightalign}
\begin{align}\label{eq:mrmr-VD}
\begin{aligned}
  \max \quad &\sum_{l\in [m]}\frac{s_l}{l^2}\\
   \text{s.t.} \quad& r = \sum_{j\in [m]}\sum_{k\in [m]} c_{jk}t_{jk}\\
  & \sum_{j\in [m]}x_j = \sum_{l\in [m]}lw_l\\
   &\sum_{l\in [m]}w_l = 1\\
   &0 \leq t_{jk}\leq x_j, \quad x_j + x_k -1 \leq t_{jk} \leq x_k && \for j, k \in [m]\\
   & s_l \leq \min\biggl\{M w_l, r+M(1-w_l)\biggr\} && \for l\in [m] \\
   & L \le \sum_{k \in [m]} x_k \le U, \quad \x, \boldsymbol{w} \in \{0, 1\}^m.
   \end{aligned}
% \tag{\textsc{mRMR-VD}}
\end{align} 
\end{tightalign}
We set $M = \vert \sum_{j,k\in [m]}c_{jk}\vert$ following \cite{mehmanchi2021solving}.

\subsection{Dataset Description}
In this section, we provide some supplementary details for the datasets used in Section~\ref{section:computation_mrmr}.

\subsubsection{Real datasets.}
\label{apx: mRMR_real_data} Here we briefly introduce several representative studies demonstrating the effectiveness of mRMR for feature selection across datasets in Section~\ref{sec: mRMR_real_data} and diverse classifiers:~\cite{singha2018adaptive} report that mRMR achieves the highest balanced average accuracy on the \texttt{Statlog} dataset when using Naive Bayes and regularized discriminant analysis. Similarly, \cite{wan2022r2ci} demonstrate that mRMR outperforms other information-theoretic feature selection methods on the \texttt{Dermatology} and \texttt{GSE28700} datasets, achieving the highest average classification accuracy with a support vector machine (SVM) classifier. In another study,
% \cmt{Jun}{change further to other words? The logic seems not ``further''.} 
\cite{gao2016variational} show that mRMR-selected features lead to low cross-validation error rates on the \texttt{Musk2} and \texttt{Lung} datasets when applying a linear SVM. Additionally, \cite{naghibi2014semidefinite} find that mRMR-based feature selection consistently yields higher average classification accuracy than joint mutual information methods on the \texttt{Lung Cancer}, \texttt{Arrhythmia}, and \texttt{CNAE-9} datasets, as measured by 10-fold cross-validation with five different classifiers. Their analysis also highlights that global optimization approaches for the mRMR problem improve both predictive performance and feature diversity. Moreover, \cite{nguyen2014effective} identify mRMR as a top-performing feature selection technique on the \texttt{Optdigits} and \texttt{Multiple Features} datasets. 

\subsubsection{Synthetic datasets.}
\label{apx: mRMR_synthetic_data}
Here we introduce the procedure to generate a synthetic dataset in Section~\ref{sec: mRMR_synthetic} as described in Section~8 of the online supplement from \cite{park2020subset}, which operates as follows:
 First, the target variable $Y = (Y_1, \ldots, Y_n)^\top$ is generated by independently sampling $Y_i \sim N(0,5)$ for $i=1,\ldots,n$. Second, the $m$ features are systematically partitioned into $m/5$ distinct groups. Within each group, an initial feature  is generated to exhibit a moderate linear relationship with $Y$, characterized by a fixed correlation coefficient $\rho=0.2$. Third, for each initially generated feature, four additional features are created to maintain strong intra-group correlations, with uniformly distributed pairwise correlation coefficients $\rho\sim \mbox{Unif}(0.5,0.8)$. Therefore, each group of five features exhibits substantial within-group correlations while maintaining moderate correlations with $Y$.

\section{Technical Proofs}
\vspace{0.5cm}
\subsection{Proof of Theorem~\ref{them:mRMR-Pers}}\label{app:proof}
Let $S_{\textsc{pers}}$ (resp. $S_{\textsc{rmc}}$) be the feasible region of ~\eqref{eqn: mRMR_full_cons} (resp.~\eqref{eq:mrmr-RMC}), and let $P_{\textsc{pers}}$ (resp. $P_{\textsc{rmc}}$) be the natural continuous relaxation of $S_{\textsc{pers}}$ (resp. $S_{\textsc{rmc}}$). Since~\eqref{eq:mrmr-RMC} is an MIP formulation of~\eqref{eq:mrmr-fp}, and it has the same objective function as~\eqref{eqn: mRMR_full_cons} does, it suffices to show that $S_{\textsc{pers}} = S_{\textsc{rmc}}$ and $P_{\textsc{pers}} \subseteq P_{\textsc{rmc}}$.

To show $P_{\textsc{pers}} \subseteq P_{\textsc{rmc}}$, we consider a point $\bar{\boldsymbol{w}} :=(\bar{\x}, \bar{\rho}, \bar{\y}, \bar{\z})$ in $P_{\textsc{pers}}$ and argue that it satisfies every linear constraints of~\eqref{eq:mrmr-RMC}. Clearly, this point satisfies all linear constraints in the first and second line of the feasible region of~\eqref{eq:mrmr-RMC}. The proof is complete by observing that the point $\bar{\boldsymbol{w}}$ also satisfies the constraints in the last line of~\eqref{eq:mrmr-RMC}, that is,
\begin{subequations}
 \label{eq:mRMR-trilinear-evelope}
\begin{align}
     0 & \leq  z_{ij} \leq y_i    \label{eq:mRMR-trilinear-evelope-1} \\ \rho^U x_j + y_i - \rho^U& \leq  z_{ij} \leq \rho^Ux_j.    \label{eq:mRMR-trilinear-evelope-2}
 \end{align}
\end{subequations}
To see that the first inequality in~\eqref{eq:mRMR-trilinear-evelope-2} is satisfied, we observe that $\bar{z}_{ij} \geq \bar{y}_{i} + \bar{y}_j - \bar{\rho} \geq \bar{y}_i + \rho^U (\bar{x}_j - 1)$, where  the second inequality holds due to $\bar{x}_j -1 \leq 0$ and the relation $\bar{y}_j \leq \rho^U \bar{x}_j $ in~\eqref{eq:mRMR_full_cons-5}. The second inequality in~\eqref{eq:mRMR-trilinear-evelope-2} is satisfied since $\bar{z}_{ij} \leq \bar{y}_j \leq \rho^U\bar{x}_j$, where the second inequality holds due to~\eqref{eq:mRMR_full_cons-5}.

Now, we can conclude that $S_{\textsc{pers}} = S_{\textsc{rmc}}$ since $P_{\textsc{pers}} \subseteq P_{\textsc{rmc}}$ implies $S_{\textsc{pers}} \subseteq S_{\textsc{rmc}}$, and, on the other hand,  $ S_{\textsc{rmc}} \subseteq P_{\textsc{pers}}$ implies $S_{\textsc{rmc}} \subseteq S_{\textsc{pers}}$. \hfill \Halmos

\subsection{Proof of Lemma~\ref{lem: review_selection_dod}}\label{apx: prof_lem1}

Let $\x = (x_1, x_2, \ldots, x_n) \in \{0,1\}^n$ be a vector of binary variables with $x_i$ modeling the selected reviews, and define
\begin{equation}
\label{eqn: def_mu}
    \mu_j(\x) :=   
 \begin{cases}
\log(P_{j}\one^\intercal \x )- \log(\d_j^\intercal \x )  &  \d_j^\intercal \x > 0\\
            \dfrac{\delta}{P_j}      & \d_j^\intercal \x = 0 
\end{cases}
\quad \for j \in [m]. 
\end{equation}
Then, it follows readily that~\eqref{eqn: Entropy_review_selection} can be expressed as 
\begin{tightalign}
    \begin{align}
    \min_{\x} \biggl\{  \sum_{j \in [m]}P_j \bigl\vert\mu_j(\x) \bigr\vert \biggm|
      \x \in \X\biggr\},
 \label{eq:IRSS_0} 
 \end{align}
 \end{tightalign}
 where $\X = \{\x \in \{0,1\}^n \mid L\leq  \sum_{i\in [n]}x_i \leq U\}$. Thus, the proof of Lemma~\ref{lem: review_selection_dod} suffices to show that~\eqref{eq:IRSS_0} is equivalent to~(\ref{eq:diff}) in the following two cases.

 \textbf{Case 1}. We consider the case where $\X_1:=\bigl\{\x \in \X \bigm |  \d_j^\intercal \x >0  \text{ for all } j\in [m] \bigr\}$ is empty. In this case, we consider $\X_2:= \{\x \in \X \mid  \one^\intercal \x = U\}$, and will prove that optimal solutions to both problems~\eqref{eq:IRSS_0} and~\eqref{eq:diff} belong to $\X_2$. Then, the proof is complete since for $\x \in \X_2$, $\mu_j(\x) = f_j(\x) - g_j(\x)$ for every $j \in [m]$. 

To prove this, denote $\mathcal{J}(\x) = \{j\mid\d_{j}^\intercal \x = 0\}$ and suppose that the optimal solution $x^*$ to  Problem~(\ref{eq:IRSS_0}) satisfies $\one^\intercal \x^*  < U$. Then we can always construct a better solution $\x^{**}$ as follows:
 \begin{equation*}
       \x^{**} = \x^* + e_{i_{j'}},
 \end{equation*}
where $i_{j'}\in S_{j'}$ and $j' \in \mathcal{J}(x^*)$. 
Now we have $\d_{j'}^\intercal \x^{**} > 0$ and the difference between the objective value of $\x^{**}$ and that of $\x^*$ is
\begin{tightalign}
    \begin{align*}
&\quad \sum_{j \in [m]}P_j \cdot \left(\bigl\vert\mu_j(\x^{**}) \bigr\vert - \bigl\vert\mu_j(\x^{*}) \bigr\vert\right)\\
& = \sum_{j \in [m]\setminus \mathcal{J}(x^*)} P_j \cdot \left(\bigl\vert\mu_j(\x^{**}) \bigr\vert - \bigl\vert\mu_j(\x^{*}) \bigr\vert\right) + \sum_{j \in \mathcal{J}(x^*) \setminus j'} P_j \cdot \left(\bigl\vert\mu_j(\x^{**}) \bigr\vert - \bigl\vert\mu_j(\x^{*}) \bigr\vert\right) \\
& \quad \quad +  P_{j'} \cdot \left(\bigl\vert\mu_{j'}(\x^{**}) \bigr\vert - \bigl\vert\mu_{j'}(\x^{*}) \bigr\vert\right)\\
& \le \sum_{j \in [m]\setminus \mathcal{J}(x^*)} P_j \bigl\vert\mu_j(\x^{**}) \bigr\vert + 0 + P_{j'} \cdot \left(\bigl\vert\mu_{j'}(\x^{**}) \bigr\vert - \bigl\vert\mu_{j'}(\x^{*}) \bigr\vert\right) \\
& = \sum_{j \in [m]\setminus \mathcal{J}(x^*)} P_j\cdot\bigl\vert\log(P_{j}\one^\intercal \x^{**} )- \log(\d_j^\intercal \x^{**} ) \bigr \vert \\
& \quad \quad + P_{j'}\cdot\left(\left\vert \log(P_{j'}\one^\intercal \x^{**} )- \log(\d_{j'}^\intercal \x^{**} ) \right\vert - \dfrac{\delta}{P_{j'}} \right)\\
& < (m+1 - \vert \mathcal{J}(x^*) \vert)\cdot\log(n) - \delta < 0,
\end{align*}
\end{tightalign}
where the second inequality holds since 
\begin{tightalign}
    \begin{align*}
\left\vert\log(P_{j}\one^\intercal \x^{**} )- \log(\d_j^\intercal \x^{**} ) \right \vert = \left\vert\log(P_{j})- \log(\dfrac{\d_j^\intercal \x^{**} }{\one^\intercal \x^{**} }) \right \vert \le \log(n) \for j \in [m]\setminus \mathcal{J}(x^{**}),
\end{align*}
\end{tightalign}
% $P_j$ and 
% $\dfrac{\d_j^\intercal \x}{\one^\intercal \x }$ 
% belong to $\left\{1, \dfrac{1}{2}, \dfrac{1}{3}, \ldots, \dfrac{1}{n}\right\}$.
and $P_{j'} < 1$. Therefore,
it contradicts the optimality of $\x^*$. Similarly, for~(\ref{eq:diff}), with the same notation, we have 
\begin{tightalign}
    \begin{align*}
&\quad \sum_{j \in [m]}P_j \cdot \left(\bigl\vert\mu_j(\x^{**}) \bigr\vert - \bigl\vert\mu_j(\x^{*}) \bigr\vert\right)\\
% & = \sum_{j \in [m]\setminus \mathcal{J}(x^*)} P_j \cdot \left(\bigl\vert\mu_j(\x^{**}) \bigr\vert - \bigl\vert\mu_j(\x^{*}) \bigr\vert\right) + \sum_{j \in \mathcal{J}(x^*) \setminus j'} P_j \cdot \left(\bigl\vert\mu_j(\x^{**}) \bigr\vert - \bigl\vert\mu_j(\x^{*}) \bigr\vert\right) \\
% & \quad +  P_{j'} \cdot \left(\bigl\vert\mu_{j'}(\x^{**}) \bigr\vert - \bigl\vert\mu_{j'}(\x^{*}) \bigr\vert\right)\\
% & \le \sum_{j \in [m]\setminus \mathcal{J}(x^*)} P_j \bigl\vert\mu_j(\x^{**}) \bigr\vert + 0 + P_{j'} \cdot \left(\bigl\vert\mu_{j'}(\x^{**}) \bigr\vert - \bigl\vert\mu_{j'}(\x^{*}) \bigr\vert\right) \\
& \le \sum_{j \in [m]\setminus \mathcal{J}(x^*)} P_j\cdot\bigl\vert\log(P_{j}\one^\intercal \x^{**} )- \log(\d_j^\intercal \x^{**} ) \bigr \vert \\
& \quad \quad + P_{j'}\cdot\left(\biggl\vert \log(P_{j'}\one^\intercal \x^{**} )- \log(\d_{j'}^\intercal \x^{**} ) \biggr\vert -  \left \vert \log(P_{j'} \one^\intercal \x^{*}) - \log(P_{j'}\cdot U) + \dfrac{\delta}{P_{j'}} \right \vert  \right)\\
% & = \sum_{j \in [m]\setminus \mathcal{J}(x^*)} P_j\cdot\bigl\vert\log(P_{j}\one^\intercal \x^{**} )- \log(\d_j^\intercal \x^{**} ) \bigr \vert \\
% & \quad + P_{j'}\cdot\left(\biggl\vert \log(P_{j'}\one^\intercal \x^{**} )- \log(\d_{j'}^\intercal \x^{**} ) \biggr\vert -  \left \vert \log(P_{j'} \one^\intercal \x^{*}) - \log(P_{j'}\cdot U) + \dfrac{\delta}{P_{j'}} \right \vert  \right)\\
& =  \sum_{j \in [m]\setminus \mathcal{J}(x^*)} P_j\cdot\bigl\vert\log(P_{j}\one^\intercal \x^{**} )- \log(\d_j^\intercal \x^{**} ) \bigr \vert \\
& \quad \quad + P_{j'}\cdot\left(\left\vert \log(P_{j'}\one^\intercal \x^{**} )- \log(\d_{j'}^\intercal \x^{**} ) \right\vert + \log(P_{j'}\cdot U) - \log(P_{j'} \one^\intercal \x^{*})  \right) -\delta \\
& < (m+2 - \vert \mathcal{J}(x^*) \vert)\cdot\log(n) - \delta \le 0,
\end{align*}
\end{tightalign}
which also contradicts the optimality of $\x^*$. Therefore, the optimal solutions to both problems~\eqref{eq:IRSS_0} and~\eqref{eq:diff} belong to $\X_2$ and problems~\eqref{eq:IRSS_0} and~\eqref{eq:diff} are equivalent in this case.
 
\textbf{Case 2}. We consider the case where $\X_1$ is not empty. Then it follows from the definitions that  for every $\x \in \X_1$, $\mu_j(\x) = f_j(\x) - g_j(\x)$ for every $j \in [m]$, and~\eqref{eq:IRSS_0} is equivalent to~\eqref{eq:diff}. Therefore, it suffices to prove that the optimal solutions to both problems~\eqref{eq:IRSS_0} and~\eqref{eq:diff} belong to $\X_1$ in this case.  
% We first prove that the optimal solutions to both problems~\eqref{eq:IRSS_0} and~\eqref{eq:diff} belong to $\X_1$. 

To prove this, suppose one solution $\x \notin \mathcal{X}_1$, then there exists some $j' \in [m]$ such that $\d_{j'}^\top \x = 0$. For problem~\eqref{eq:IRSS_0}, the objective 
\[
\sum_{j \in [m]}P_j \bigl\vert\mu_j(\x) \bigr\vert \ge P_{j'} \cdot \dfrac{\delta}{P_{j'}} \ge \delta.
\]
% It suffices to find a solution $\x \notin \mathcal{X}_1$ whose objective is lower than $n\log(n)$.
However, for any solution $\x \in \mathcal{X}_1$, the objective of Problem~\eqref{eq:IRSS_0}
\begin{tightalign}
    \begin{align*}
\sum_{j \in [m]}P_j \bigl\vert\mu_j(\x) \bigr\vert & = \sum_{j \in [m]}P_j \bigl\vert\log(P_{j}\one^\intercal \x )- \log(\d_j^\intercal \x )  \bigr\vert  \\
& = \sum_{j \in [m]}P_j \left\vert\log(P_{j}) - \log(\dfrac{\d_j^\intercal \x}{\one^\intercal \x }) \right\vert \\
& \le \sum_{j \in [m]}P_j \log(n)  \le m\log(n) < \delta.
\end{align*}
\end{tightalign}
Therefore, the optimal solutions to problems~\eqref{eq:IRSS_0} fall into $\mathcal{X}_1$. Similarly, for~\eqref{eq:diff}, if $\x \notin \mathcal{X}_1$ and with the same notation, we have
\begin{tightalign}
    \begin{align*}
 \sum_{j \in [m]}P_j \bigl\vert f_j(\x)- g_j(\x) \bigr\vert 
& \ge P_{j'}\cdot \biggl\vert \log(P_{j'}\one^\intercal \x ) - \log(P_{j'}\cdot U) + \frac{\delta}{P_{j'}}  \biggr\vert \\
& = \delta -  P_{j'}\cdot \bigl\vert \log(P_{j'}\one^\intercal \x ) - \log(P_{j'}\cdot U)  \bigr\vert.
\end{align*}
\end{tightalign}
In contrast, for any solution $\x \in \mathcal{X}_1$, the objective of~\eqref{eq:diff} is same to~\eqref{eq:IRSS_0}, which is less than $\delta -  P_{j'}\cdot \bigl\vert \log(P_{j'}\one^\intercal \x ) - \log(P_{j'}\cdot U)  \bigr\vert$. Therefore, the optimal solutions to~\eqref{eq:IRSS_0} also fall into $\mathcal{X}_1$ and problems~\eqref{eq:diff} and~\eqref{eq:IRSS_0} are equivalent in this case. \hfill \Halmos

\subsection{Proof of Proposition~\ref{prop: review_selection_submodular}} 
Due to Theorem~3.3 in~\cite{tawarmalani2013explicit}, it suffices to that both $f_j(\cdot)$ and $g_j(\cdot)$ are submodular functions. To establish this, we will use the following lemma. 
\begin{lemma}\label{lemma:submodular}
 Consider a composite function $h:\{0,1\}^n \setminus \{\boldsymbol{0}\} \to \R$ defined as $h(\x) = \varphi (\boldsymbol{\alpha}^\intercal \x)$, where $\varphi(\cdot)$ is a concave function defined over  the positive numbers and $\boldsymbol{\alpha}$ is a vector of positive numbers. Let $\{a_1, a_2, \ldots, a_N \}$ denote the range $\{\boldsymbol{\alpha}^\intercal \x \mid \x \in \{0,1\}^n \setminus \{\boldsymbol{0}\}\}$ such that $0 < a_1< a_2 < \cdots < a_N$, and let 
 \begin{equation}\label{eq:lemma-condition-1}
     \tau:= \varphi(a_1) - \frac{\varphi(a_2) - \varphi(a_1)}{a_2-a_1} a_1. 
 \end{equation}
 Then, an extension $\bar{h}:\{0,1\}^n \to \R$ of $h(\cdot)$ is submodular if $\bar{h}(\boldsymbol{0}) \leq \tau$. 
\end{lemma}
\noindent {\bf Proof.} 
Let $a_0 = 0$, and let $\bar{\varphi}:\{a_0, a_1, \ldots, a_N\} \to \R$ be a function such that $\bar{\varphi}(y) = \varphi(y)$  if $y \in \{a_1, a_2, \ldots, a_N\}$ and $\bar{\varphi}(a_0) \leq \tau$. It follows readily that $\bar{h}(\x) = \bar{\varphi}(\x)$ for every $\x \in \{0,1\}^n$. Moreover, for $y', y'' \in \{a_0, a_1, \ldots, a_N\}$ with $y'\geq y''$ and $\delta$ such that $y'+ \delta$ and $y'' + \delta$ belong to $\{a_0, a_1, \ldots, a_N \}$, we have
\begin{equation}\label{eq:EC-prop1-1}
    \bar{\varphi}(y'+ \delta) -\bar{\varphi}(y') \leq \bar{\varphi}(y''+ \delta) -\bar{\varphi}(y''). 
\end{equation}
% The inequality holds as follows. The definition of $\bar{\varphi}(\cdot)$ and the concavity of $\varphi(\cdot)$ imply that 
% \[
% \frac{\bar{\varphi}(a_i) - \bar{\varphi}(a_{i+1})}{a_i-a_{i+1}} \geq \frac{\bar{\varphi}(a_j) - \bar{\varphi}(a_{j+1})}{a_j-a_{j+1}} \qquad  \for \, i \leq j. 
% \]
% This implies that 
Let $\x'$ and $\x''$ be two vectors in $\{0,1\}^n$. Let $y':=\boldsymbol{\alpha}^\intercal\x'$,  $y'':=\boldsymbol{\alpha}^\intercal( \x' \wedge \x'')$ and $\delta = \sum_{i \in [n]} \alpha_i \max\{0, x''_i -x'_i \}$. Then, we have 
\[
\begin{aligned}
\bar{h}(\x'\vee \x'') + \bar{h}(\x'\wedge \x'')  -  \bar{h}(\x') - \bar{h}(\x'') &=  \bigl(\bar{h}(\x'\vee \x'') -  \bar{h}(\x')\bigr) - \bigl(  \bar{h}(\x'') - \bar{h}(\x'\wedge \x'') \bigr)  \\    
& =\bigl( \bar{\varphi}(y' + \delta) - \bar{\varphi}(y')\bigr) - \bigl( \bar{\varphi}(y'' + \delta) - \bar{\varphi}(y'')\bigr) \\
& \leq 0,
\end{aligned}
\]
where the second equality holds by definition, and the inequality follows from~\eqref{eq:EC-prop1-1}. This shows the submodularity of $\bar{h}(\cdot)$.  \hfill \Halmos

By Lemma~\ref{lemma:submodular}, the submodularity of $f_j(\cdot)$ holds since the condition~\eqref{eq:lemma-condition-1} is satisfied:  
\[
f_j(\boldsymbol{0}) = 2 \cdot \phi_j(1) - \phi_j(2) = \phi_j - \frac{\phi_j(2) - \phi_j(1)}{2-1} 1.
\]
The submodularity $g_j(\cdot)$ holds since condition~\eqref{eq:lemma-condition-1} is satisfied:
\[
g_j(\boldsymbol{0}) =\log(P_j U) - \dfrac{n\log(n)}{P_j} \le \log(P_j U)- \dfrac{U\log(U)}{P_j} \le -\log(2) = \psi_j(1) - \dfrac{\psi_j(2) - \psi_j(1)}{2-1} \cdot 1,
\]
where two inequalities hold when $U \geq 2$. \hfill \Halmos

%  Since $\psi_j(\cdot)$ is concave over the positive numbers and the vector $\one_{\vert S_j \vert \times 1}$ is also positive, by Lemma~\ref{lemma:submodular}, the proof for the submodularity of $g_j(\cdot)$ suffices to show that 
% \[
% g_j(\boldsymbol{0}) =\log(P_j U) - \dfrac{n\log(n)}{P_j} \le \psi_j(1) - \dfrac{\psi_j(2) - \psi_j(1)}{2-1} \cdot 1   = -\log(2).
% \]
% As 
% \[\log(P_j U) - \dfrac{n\log(n)}{P_j}\le \log(P_j U)- \dfrac{U\log(U)}{P_j} \le -\log(2)\]
% always holds when $U\geq 2$, the proof is complete. 
% which is reasonable in practice. 
% We finally obtain 
 % \begin{equation*}
    % \psi_j(0) =  \phi_j(U) - \dfrac{n\log(n)}{P_j}.
% \end{equation*}

\subsection{Proof of Proposition~\ref{prop:conc}}
To prove this result, we need the following lemma. 
\begin{lemma}
\label{lem: conc_env}
    Let $\boldsymbol{\alpha}$ be a vector in $\{0,1\}^n$ and define $S:= \bigl\{i \in [n] \bigm | \alpha_i \neq 0\bigr\}$. Consider a composite function $h:\{0,1\}^n \to \R$ defined as $h(\x) = \psi(\boldsymbol{\alpha}^\intercal \x)$, where $\psi : \{0,1,\ldots, \vert S \vert\}  \to \R$. If for $k \in [\vert S \vert]$,
    \begin{equation}\label{eq:lemma-condition-2}
        \psi(y) \leq \bigl(\psi(k) - \psi(k-1)\bigr) y + k \psi(k-1) - (k-1)\psi(k) \for y \in \{0,1, \ldots, \vert S\vert\},
    \end{equation}
    then
    \[
    \conc(h)(\x) = \min_{k \in [\vert S \vert]} \biggl\{  \underbrace{\left[\psi(k) - \psi(k-1)\right]\sum_{i \in S} x_i + k \psi(k-1)-(k-1)\psi(k)}_{=:\ell_k(\x)}  \biggr\} \quad \for \, \x \in [0,1]^n.
    \]
\end{lemma}
\noindent {\bf Proof.}  Let $L(\x):= \min_{k \in [\vert S \vert]} \{\ell_k(\x)\}$ for $\x \in [0,1]^n$. First, we show that $L(\x) \geq \conc(h)(\x)$ for every $\x \in [0,1]^n$. Let $y:=\sum_{i \in S}x_i$, and then, for $k \in [\vert S \vert ]$, 
\[
\ell_k(\x) = \left[\psi(k) - \psi(k-1)\right] y + k \psi(k-1) - (k-1)\psi(k) \geq \psi(y) = h(\x) \quad \for \x \in \{0,1\}^n,
\]
where the first and last equality hold by definitions, and the inequality holds by the assumption on $\psi(\cdot)$.  Therefore, $\ell_k(\x) \geq \conc(h)(\x)$ for $\x \in [0,1]^n$ since $\ell_k(\cdot)$ is a concave over-estimator of $h(\cdot)$ and $\conc(h)(\cdot)$ is the tightest concave over-estimator of $h(\cdot)$. Hence, $L(\x) \geq \conc(h)(\x)$.

Next, we show that $L(\x) \leq \conc(h)(\x)$ for every $\x \in [0,1]^n$. 
For $ k \in [\vert S \vert]$, consider a polytope $P_k := \bigl\{\x \in [0,1]^n \bigm| k-1\leq  \sum_{i \in S} x_i  \leq k \bigr\}$. It can be shown that the set of vertices of $P_k$, denoted as $\text{vert}(P_k)$ is $V_k-1 \cup V_k$, where for $k \in \{0,1, \ldots, \vert S \vert \}$, 
\[
V_k :=\biggl\{\x \in \{0,1\}^n \biggm| \sum_{i \in S} x_i =k\biggr\}. 
\]
For every $k \in [\vert S\vert]$, $\ell_k(\x) = h(\x)$ for every $x \in \text{vert}(P_k)$ since for every $\x \in V_{k-1}$,
\[
\ell_k(\x) = (\psi(k) - \psi(k-1))(k-1) +k\psi(k-1) - (k-1)\psi(k) = \psi(k-1) = h(\x),
\]
and for every $\x \in V_{k}$,
\[
\ell_k(\x) = (\psi(k) - \psi(k-1))k +k\psi(k-1) - (k-1)\psi(k) = \psi(k) = h(\x). 
\]
Now, let $\x'$ be a point in $[0,1]^n$. Then, there exists $k'$ such that $\x' \in P_{k'}$, and thus there exists a convex multiplier $\boldsymbol{\lambda}$ such that $\x' = \sum_{v \in \text{vert}(P_{k'})}v \lambda_v $. It turns out that  
\[
L(\x') \leq \ell_{k'}(\x')  = \sum_{v \in \text{vert}(P_{k'})}\lambda_v\ell_{k'}(v) = \sum_{v \in \text{vert}(P_{k'})}\lambda_v h(v) \leq \conc(h)(\x'),
\]
where the first inequality holds since $\ell_{k'}(\cdot)$ is one of affine functions defining $L(\cdot)$, the first inequality holds due to the linearity of $\ell_{k'}(\cdot)$, the second equality follows from the above discussion, and the last inequality follows from the definition of the concave envelope. \hfill \Halmos

By Lemma~\ref{lem: conc_env}, Proposition~\ref{prop:conc} holds when ~\eqref{eq:lemma-condition-2} is satisfied by $\phi_j(\cdot)$ and $\psi_j(\cdot)$. For $\phi_j(\cdot)$, to show this, we first prove that, 
\[
\phi_j(k + 1) - \phi_j(k) \le \phi_j(k) - \phi_j(k - 1) \for k \in [n - 1].
\]
It holds for $k \in [n -1]\setminus \{1\}$ since $\phi_j(\cdot)$ is concave over the domain of positive numbers. For $k = 1$, we have $\phi_j(k + 1) - \phi_j(k) = \log(2 P_j) - \log(P_j) = \phi_j(k) - \phi_j(k-1)  $ also holds by definition of  $\phi_j(\cdot)$.

Then when $y \ge k$,  we have 
\begin{tightalign}
\begin{align*} 
\phi_j(y) - \phi_j(k -1) & =  \sum_{i = k}^{y} \left( \phi_j(i) - \phi_j(i - 1) \right) \\
& \le  \sum_{i = k}^{y} \left( \phi_j(k) - \phi_j(k - 1) \right) \\
&= (y - k + 1)\left( \phi_j(k) - \phi_j(k - 1) \right),
\end{align*} 
\end{tightalign}
which is equivalent to~\eqref{eq:lemma-condition-2}. When $y \le k - 1$, we have 
\begin{tightalign}
\begin{align*} 
\phi_j(k) - \phi_j(y) & =  \sum_{i = y + 1}^{k} \left( \phi_j(i) - \phi_j(i - 1) \right) \\
& \ge  \sum_{i = y + 1}^{k} \left( \phi_j(k) - \phi_j(k - 1) \right) \\
&= (k - y)\left( \phi_j(k) - \phi_j(k - 1) \right),
\end{align*} 
\end{tightalign}
which is also equivalent to~\eqref{eq:lemma-condition-2}. Therefore,~\eqref{eq:lemma-condition-2} always holds by $\phi_j(\cdot)$. Similar proof procedure can be applied for $\psi_j(\cdot)$.  \hfill \Halmos

\subsection{Proof of Theorem~\ref{theorem:envelope}}
By Lemma~\ref{lem: review_selection_dod}, we obtain that~\eqref{eqn: Entropy_review_selection}  is equivalent to~\eqref{eq:diff}. For each $j \in [m]$, let $F_j$ (resp. $G_j$) be the graph of $f_j(\cdot)$ (resp. $g_j(\cdot)$), that is,
\[
F_j = \bigl\{(\x, s_j) \bigm| s_j = f_j(\x),\ \x \in \{0,1\}^n \bigr\} \quad \text{ and } \quad G_j = \bigl\{(\x, t_j) \bigm| t_j = g_j(\x),\ \x \in \{0,1\}^n\bigr\} .
\]
Since $\{0,1\}^n$ is the set of vertices of $[0,1]^n$, $\conv(f_j)(\x) = f_j(\x) = \conc(f_j)(\x)$  and $\conv(g_j)(\x) = g_j(\x) = \conc(g_j)(\x)$ for every $x\in \{0,1\}^n$. In  other words, we have 
\[
\begin{aligned}
F_j = \bigl\{(\x, s_j) \bigm| \conc(f_j)(\x) \geq  s_j \geq \conv(f_j)(\x),\ \x \in \{0,1\}^n \bigr\} \\
G_j = \bigl\{(\x, t_j) \bigm| \conc(g_j)(\x) \geq  t_j \geq \conv(g_j)(\x),\ \x \in \{0,1\}^n \bigr\}. 
\end{aligned}
\]
Now, using  the explicit descriptions of envelopes in Propositions~\ref{prop: review_selection_submodular} and ~\ref{prop:conc}, and the fact that  the constraint $\mu_j \geq \vert s_j - t_j\vert$ is equivalent to $\mu_j \geq s_j - t_j$  and $\mu_j \geq t_j -s_j$, we can conclude that~\eqref{eq:dataformulation} is an MIP formulation~\eqref{eqn: Entropy_review_selection}. \hfill \Halmos

\end{document}